\documentclass[preprint,12pt]{elsarticle}

\usepackage[letterpaper,top=1.75cm,bottom=1.75cm,left=1.75cm,right=1.75cm,marginparwidth=1.75cm]{geometry}

\usepackage{bm}
\usepackage{amsmath}
\usepackage{amsfonts}
\usepackage{graphicx}[realboxes]
\usepackage[colorlinks=true, allcolors=blue]{hyperref}
\usepackage{placeins}
\usepackage{hyperref}
\usepackage{multirow}
\usepackage{caption}
\usepackage{subcaption}
\usepackage{booktabs}
\usepackage{siunitx}
\usepackage{rotating}
\hypersetup{pdfauthor=author}
\usepackage{adjustbox}

\usepackage[ruled,vlined,linesnumbered]{algorithm2e}
\usepackage{algpseudocode}

\SetCommentSty{mycommfont}

\newcommand{\R}{\mathbb{R}}
\newcommand{\bA}{\bm{A}} 
\newcommand{\bE}{\bm{E}} 
\newcommand{\bS}{\bm{S}} 

\newcommand{\bOmega}{\bm{\Omega}} 

\DeclareMathOperator*{\argmin}{arg\,min}
\DeclareMathOperator{\tr}{tr}
\DeclareMathOperator{\spantext}{span}

\newcommand{\pinv}{\dagger} 

\newtheorem{theorem}{Theorem}
\newtheorem{lemma}{Lemma}

\journal{Elsevier}

\begin{document}

\begin{frontmatter}

        \title{Online randomized interpolative decomposition with {\it a posteriori} error estimator for temporal PDE data reduction}
        \author[1]{Angran Li}
        \author[2]{Stephen Becker}
        \author[1]{Alireza Doostan \corref{cor}}

        \address[1]{Smead Aerospace Engineering Sciences, University of Colorado Boulder, 80309, Boulder, USA}
        \address[2]{Applied Mathematics, University of Colorado Boulder, 80309, Boulder, USA}

        \cortext[cor]{Corresponding author. E-mail: alireza.doostan@colorado.edu.}

        \begin{abstract}
            
            Traditional low-rank approximation is a powerful tool to 
            compress the huge data matrices that arise in simulations of partial differential equations (PDE), but suffers  from high computational cost and requires several passes over the PDE data. The compressed data may also lack interpretability thus making it difficult to identify feature patterns from the original data. To address these issues, we present an online randomized algorithm to compute the interpolative decomposition (ID) of large-scale data matrices {\em in situ}. Compared to previous randomized IDs that used the QR decomposition to determine the column basis, we adopt a streaming ridge leverage score-based column subset selection algorithm that dynamically selects proper basis columns from the data and thus avoids an extra pass over the data to compute the coefficient matrix of the ID. In particular, we adopt a single-pass error estimator based on the non-adaptive Hutch++ algorithm to provide real-time error approximation for determining the best coefficients. As a result, our approach only needs a single pass over the original data and thus is suitable for large and high-dimensional matrices stored outside of core memory or generated in PDE simulations. A strategy to improve the accuracy of the reconstructed data gradient, when desired, within the ID framework is also presented. We  provide numerical experiments 
            on turbulent channel flow and ignition simulations, and on the NSTX Gas Puff Image dataset, comparing our algorithm with the offline ID algorithm to demonstrate its utility in real-world applications.
        \end{abstract}


    
\end{frontmatter}

\section{Introduction}
Recent advances in parallel computing have led to an increasing interest in addressing challenges in industrial application, scientific design, and discovery using computer simulations. In particular, simulations based on partial differential equations (PDEs) have grown in size and complexity, and generate a considerable amount of data at fast rates. However, the speed of the current input/output (I/O) filesystem does not keep up with the data generation speed thus limiting the improvement of computation time. It is also challenging to store the entire simulation data due to the storage size limitation of current filesystems. Therefore, it is essential to develop compression algorithms that can efficiently reduce the size of large PDE data streams for storage or visualization.

Our setup is that $\bA = [\bm{a}_1, ..., \bm{a}_n] $ is a $m \times n$ matrix that contains $n$ snapshots of the PDE solution, each of size $m$,
i.e., $\bm{a}_t$ is the discretized solution at some discrete time step $t$. We are most interested in the case when $\bA$ is so large that it is impossible to store all at once in the main memory of a computer, and we receive incremental updates of new columns after discarding old columns, as in a time-dependent PDE simulation. Our method is applicable for both $m\ge n$ and $m<n$ cases.


Traditionally, low-rank approximation has served as a robust and powerful tool for data compression tasks.
The goal of low-rank approximation is to find a rank $k$ matrix that is as close as possible to $\bA$ with distance typically measured using the spectral or Frobenius norms. The Eckart-Young theorem \cite{eckart1936approximation} states that the best rank-$k$ approximation of $\bA$ in the spectral or Frobenius norm can be provided via the singular value decomposition (SVD):
\begin{equation}
        \bA \approx \bm{U}_k \bm{\Sigma}_k \bm{V}_k^T,
\end{equation}
where $\bm{\Sigma}_k$ is a diagonal matrix whose entries are the $k$ largest singular values of $\bA$, and $\bm{U}_k$ and $\bm{V}_k$ contain the first $k$ left and right singular vectors, respectively. Instead of storing the $mn$ numbers of $\bA$, the rank $k$ matrices $\bm{U}_k, \bm{\Sigma}_k, \bm{V}_k$ can be stored with only $k(m+n+1)$ numbers and thus achieves compression of $\bA$ if $k$ is small enough; variants of this approach that also further compress $\bm{U}_k, \bm{\Sigma}_k, \bm{V}_k$ have been proposed in the literature, e.g.,
\cite{LowrankCompressio_Pilanci2023}. 
However, it is challenging for the SVD to tackle large matrices due to its expensive computational cost in both time and memory. Although the high cost can be reduced using iterative algorithms, such as the power method or Krylov methods, these methods require access to the matrix multiple times
so they are not suitable for a streaming scenario.

The interpolative decomposition (ID) is another low-rank approximation approach that uses a sub-matrix of $\bA$ (a subset of its rows or columns) as the basis to approximate the original data. Though the ID is generally less accurate than the SVD \cite{voronin2017efficient}, the decomposition results are more interpretable since it utilizes a subset of $\bA$ as basis vectors to approximate the original matrix. For matrices with important properties such as sparsity or non-negativity, one can specifically select bases that preserve these properties and while still obtaining a compact representation of the original data. The ID method has been applied in compressing computational fluid dynamics data \cite{dunton2020pass,pacella2022task} and computational electromagnetic problems \cite{pan2012fast, huang2016efficient}, as well as in multi-fidelity uncertainty quantification~\cite{hampton2018practical,skinner2019reduced,fairbanks2020bi}. 
The idea of column ID is to approximate the matrix $\bA$ as a product of a subset of its columns (also called a column skeleton) $\bA(:,\mathcal{J})$, with $\mathcal{J}$ containing a subset of $k$ indices from $\{1,2,...,n\}$, and a coefficient matrix $\bm{P}\in \R^{k \times n}$ which reconstructs the matrix via
\begin{equation}
\label{eqn:id}
        \bA \approx \bA(:,\mathcal{J})\bm{P},
\end{equation}
where we use MATLAB notation to specify submatrices. 
A standard way to obtain the ID of $\bA$ is to first compute the column-pivoted QR (CPQR) factorization of $\bA$ and use the top $k$ pivot indices to obtain $\mathcal{J}$, and then compute $\bm{P}$ by minimizing the Frobenius norm $\|\bA - \bA(:,\mathcal{J})\bm{P}\|^2_F$. However, the CPQR factorization requires multiple passes over $\bA$, and thus ID shares a similar issue as SVD in handling large-size matrices. 

Randomized methods based on Johnson-Lindenstrauss random projection \cite{liberty2007randomized, halko2011finding, yu2017single} have garnered significant attention for their capability to handle the low-rank approximation of large matrices with reduced cost and fewer passes. 
They can also output a low-rank approximation with a (1 + $\epsilon$) factor of the optimal error when measured using the Frobenius norm and thus have strong theoretical guarantees to ensure that the approximation error is controlled, making them suitable for applications where accuracy is critical. In this study, we are specifically interested in the randomized ID method \cite{liberty2007randomized}. The key idea is first to use random projection to create a sketch matrix $\bS \in \mathbb{R}^{\ell \times n} (\ell \ll m)$ as
\begin{equation}
        \label{eq:random_sketch}
        \bS = \bOmega \bA = [\bOmega\bm{a}_1, ..., \bOmega\bm{a}_n], 
\end{equation}
where $\bOmega \in \mathbb{R}^{\ell \times m}$ is the random projection matrix and $\bS$ can be treated as a spanning set for the range of a matrix $\bA$.
If $\ell$ is sufficiently small, then $\bS$ \emph{can} be stored all at once in main memory, unlike $\bA$, and 
since $\bOmega$ can be applied to each column of $\bA$, the random projection can be computed on streaming data.
Then, CPQR is applied to $\bS$ to obtain the $k$ column indices in $\mathcal{J}$ as well as the column skeleton $\bA(:, \mathcal{J})$. Finally, the coefficient matrix $\bm{P}$ can be computed by minimizing $\|\bA-\bA(:, \mathcal{J})\bm{P}\|_F^2$.
These last two steps can be combined in a single pass, but this is still in addition to the first pass to compute $\bS$, so overall this is 
a two-pass method. Compared to traditional, deterministic methods, the randomized ID algorithm is pass-efficient and it also reduces the time and space complexity for computing the coefficient matrix since the CPQR is applied to the smaller sketch matrix and involves fewer floating-point operations (\textit{flops}). However, the difference between two passes and one pass is great: the 
two-pass randomized ID method still needs 
either
an extra simulation run to obtain the column basis
or saving data to secondary or tertiary memory. 
Hence, our motivation in this paper is to develop an online method that collects the column skeleton and computes the coefficient matrix in a single pass. 

We note that a closely related method, column subset selection (CSS), 
has been explored in previous work to identify a subset of columns that serves as the optimal basis to depict the original matrix~\cite{boutsidis2009improved,boutsidis2014near,deshpande2010efficient, paul2015column, tropp2009column, mahoney2011randomized, civril2012column}. In more detail, given a matrix $\bA$, CSS aims to find a small subset $\bm{C}$ of its columns that minimizes the projection error of the matrix to the span of the selected columns:
\begin{equation}
        \underset{\bm{X}}{\min} \|\bm{C} \bm{X}- \bA\|_F^2 =  \|\bm{C} \bm{C}^{\pinv}\bA- \bA\|_F^2,
\end{equation}
where $\bm{C}^{\pinv}$ denotes the pseudoinverse of $\bm{C}$.
The CSS and ID problems are similar since $\bm{C}$ comes from $\bA$ and one possible ID solution can be obtained by setting $\bA(:,\mathcal{J}) = \bm{C}$ and $\bm{P} = \bm{X}$. In \cite{cohen2017input, bhaskara2019residual}, the authors developed online CSS algorithms using residual norm and ridge leverage score-based random sampling which only require one pass through the data.
However, most online CSS methods only provide the column subset but still need an extra pass over the entire data to build the coefficient matrix for reconstruction.
In this regard, the randomized and online CSS algorithms complement each other well since the randomized algorithm can efficiently compute the coefficients using the sketch matrix given the selected column basis from the online CSS. This inspires us to combine both algorithms in a single framework and develop an online algorithm to handle the streaming data and obtain the ID matrices in one PDE simulation run. 



%
%
\subsection{Contributions}

In this study, we focus on addressing data compression in large-scale PDE simulations where the PDE solution at each time step arrives in a streaming fashion. 
The main contribution of this paper is the development of a novel online randomized column ID approach to perform temporal compression of large-scale simulation data. In particular, we make the following contributions:
\begin{itemize}
        \item Extending the ridge leverage score-based CSS method to compute the column basis in one pass~\cite{cohen2017input}, we use a random projection sketching technique to simultaneously compute the coefficient matrix without any additional passes through the data;
        \item We consider four different approaches to compute the coefficients 
        and adaptively choose the best one by evaluating the approximation error of different coefficients. Specifically, a single-pass error estimator based on the non-adaptive Hutch++ (NA-Hutch++) algorithm \cite{meyer2021hutch++, jiang2021optimal} is adopted to provide real-time error approximation for determining the best coefficients; 
        \item To accommodate datasets where the reconstruction accuracy of the gradient field is important, we also derive a fast gradient estimation and use the gradient information in both CSS and coefficient computation processes;
        \item Demonstrate the effectiveness of our approach by evaluating the compression performance of time-dependent scientific data obtained from several applications, namely 
        the direct numerical simulation of a turbulent channel flow~\cite{perlman2007data,li2008public,graham2016web}, an ignition simulation, and the NSTX Gass Puff Image (GPI) data~\cite{choi2020data, zhao2020sdrbench, SDRBench_git}.
\end{itemize}

The remainder of the manuscript is organized as follows. In Section \ref{sec:Background}, we provide appropriate background and context for the randomized ID and review two online CSS algorithms using residual-based and ridge leverage score-based sampling, respectively. In Section \ref{sec:Method}, we outline the online randomized ID algorithm proposed in this work. In Section \ref{sec:Result}, we present the numerical experiments for the proposed method applied to three different simulation datasets. Finally, in Section \ref{sec:Conclusion}, we summarize the findings of this work and suggest future research.

\section{Background}
\label{sec:Background}

In this section, we review the standard CPQR-based ID and online residual-based CSS algorithms that are used to build our algorithm. We use the following notation. Matrices and vectors are denoted by capital and lowercase bold letters, respectively. The $j^\text{th}$ column of $\bA$ is expressed as $\bm{a}_j$ and the $(i, j)^\text{th}$ entry of $\bA$ as $a_{ij}$.
For a matrix $\bA$, using MATLAB notation, the submatrix $\bA(:, \mathcal{J}) = \bA_{\mathcal{J}}$ represents those columns from $\bA$ with indices given by the set $\mathcal{J}$. $\bA^{\pinv}$ denotes the Moore-Penrose inverse of $\bA$. The symbol $\bA_k$ denotes the best rank-$k$ approximation of $\bA$ with respect to the Frobenius norm. In general, $\bm{C}$ denotes the selected column subset, equivalent to $\bA_{\mathcal{J}}$, while $\bOmega$ denotes the random projection matrix, and $\bS$ denotes the matrix sketch of $\bA$, $\bS=\bOmega \bA$.

\subsection{ID based on CPQR}
We first review the computation of low-rank column ID via CPQR factorization, which is later used as the benchmark to evaluate the performance of our algorithm.
The column pivoted (thin/reduced) QR factorization of $\bA\in \R^{m \times n}$ (assuming $m\ge n$) returns
\begin{equation}
        \bA \bm{Z} = \bm{QR},
\end{equation}
where $\bm{Z}$ is a permutation matrix, $\bm{Q} \in \R^{m \times n}$ has orthogonal columns and $\bm{R} \in \R^{n \times n}$ is an upper triangular matrix. $\bm{Q}$ and $\bm{R}$ can be partitioned into blocks such that
\begin{align}
        \bm{Q} & = \begin{bmatrix} \bm{Q}_k & \bm{Q}_{n-k} \end{bmatrix},                                 \\
        \bm{R} & = \begin{bmatrix} \bm{R}_{11} & \bm{R}_{12}\\ \bm{0} & \bm{R}_{22}\end{bmatrix}, 
\end{align}
where $k < \min(m,n)$, $\bm{Q}_k \in \R^{m \times k}$ contains the first $k$ orthogonal columns of $\bA$, $\bm{Q}_{n-k} \in \R^{m \times (n-k)}$, $\bm{R}_{11} \in \R^{k \times k}$, $\bm{R}_{12} \in \R^{k \times (n-k)}$, and $\bm{R}_{22} \in \R^{(n-k) \times (n-k)}$. The rank $k$ approximation of $\bA$ is obtained from the reduced CPQR as
\begin{align}
        \bA \bm{Z} & \approx \bm{Q}_k \begin{bmatrix} \bm{R}_{11} & \bm{R}_{12} \end{bmatrix} \\ &= \bm{Q}_k\bm{R}_{11} \begin{bmatrix} \bm{I}_{k} & \bm{R}_{11}^{-1} \bm{R}_{12} \end{bmatrix}.
\end{align}
Letting the index set $\mathcal{J}$ point to the first $k$ pivoted vectors of $\bA$, we have
\begin{align}
       \bA_\mathcal{J} &= \bm{Q}_k\bm{R}_{11}, \\
        \bA &\approx \bA_{\mathcal{J}} \begin{bmatrix} \bm{I}_{k} & \bm{R}_{11}^{-1} \bm{R}_{12} \end{bmatrix} \bm{Z}^T.
\end{align}
The column ID of $\bA \approx \bA_{\mathcal{J}} \bm{P}$ is obtained by setting $\bm{P} = \begin{bmatrix} \bm{I}_{k} & \bm{R}_{11}^{-1} \bm{R}_{12} \end{bmatrix} \bm{Z}^T$, which is equivalent to computing $\bm{P}$ by minimizing $\|\bA-\bA(:, \mathcal{J})\bm{P}\|_F^2$. 

\subsection{Review of online CSS algorithms}
\label{subsec:online_CSS_review}
To develop an online ID method, it is crucial to simultaneously determine the column basis and compute the current coefficient matrix in a streaming setting. Both \textit{online} and \textit{streaming} CSS methods are promising since they each make decisions on the column as it arrives, leading to the selected columns being immediately used to compute a coefficient matrix without another pass through the data. Therefore, we consider \textit{online} and \textit{streaming} to be equivalent terms and use them interchangeably throughout this paper. In the following sections, we provide a review of two online CSS algorithms using residual-based and ridge leverage score-based sampling, respectively.
\subsubsection{Online CSS using residual-based sampling}
\label{subsubsec:residual_CSS}
\newcommand{\pre}{\text{pre}}
\newcommand{\current}{\text{current}}
We first review an online CSS algorithm using residual-based sampling proposed in~\cite{bhaskara2019residual}. The basic idea is to maintain a subspace corresponding to the span of previously selected columns while collecting new columns based on their residual norm by projecting to the subspace. The arriving sequence of columns is partitioned into steps where we update the column subset until the residual norm based on the current column subset reaches a threshold.
To describe the method, we first denote the previously and currently selected column subsets as $\bm{C}_{\pre}$ and $\bm{C}_{\current}$, respectively. $\Pi^{\perp}_{\bm{C}}$ denotes the orthogonal projection onto the subspace orthogonal to the column space of $\bm{C}$, and $\sigma$ denotes the ``residual mass'' that measures the accumulated normalized residual of approximating new data columns with the currently selected columns in each CSS step. The algorithm takes the data matrix $\bA$, an estimate of the target rank $k$, and a parameter satisfying $\xi \geq \|\bA-\bA_{k}\|^2_F$ as input. Both $\bm{C}_{\pre}$ and $\bm{C}_{\current}$ are initialized to the empty set, and $\sigma$ is set to 0. When a new data column $\bm{a}$ arrives, we first compute its probability to be selected into the subset $\bm{C}_{\current}$ as 
\begin{equation}
    p_a := \dfrac{k\|\Pi^{\perp}_{\bm{C}_{\pre}}\bm{a}\|_2^2}{\xi}.
\end{equation}
This new column is then added to $\bm{C}_{\current}$ with probability
$\text{min}(p_a,1)$. 
The probability $p_a$ is also added to $\sigma$ as a collection of residual mass in this step. The process continues until $p_a \geq 1$ or $\sigma \geq 1$, indicating that the residual norm of using the current subset to approximate the new observed columns reaches the error threshold. Therefore, the previously selected columns $\bm{C}_{\pre}$ are not sufficient to provide an accurate approximation of the new columns observed in this step, and thus an update with the newly selected columns $\bm{C}_{\current}$ is performed:
\begin{equation}
    \bm{C}_{\pre} \gets \bm{C}_{\pre} \cup \bm{C}_{\current}.
\end{equation}
Note that $\bm{C}_{\pre}$ is not updated unless it has built up enough residual mass $\sigma$, and thus $\sigma$ can be used to bound the number of selected columns \cite{bhaskara2019residual}. At this point, $\sigma=0$ and $\bm{C}_{\current} = \emptyset$ are reset to start a new step of collection columns. The process continues until the entire data matrix is loaded to obtain the CSS of $\bA$ with selected columns $\bm{C} = \bm{C}_{\pre}$. 
As proved in \cite{bhaskara2019residual}, given an integer $k \ge 1$ and a fixed parameter $\xi \ge \|\bA - \bA_k\|^2_F$, with probability $\geq 1- \delta$, the algorithm provides a subset $\bm{C}$ such that
\begin{equation}
\|\Pi^{\perp}_{\bm{C}}\bA\|_F^2 \leq \xi \cdot \mathcal{O}\left(\log \dfrac{\|\bA\|^2_F}{\xi} + \dfrac{\log(1/\delta)}{k}\right).
\end{equation}

The selected columns can then be used as the basis to compute the coefficient matrix $\bm{P}$ and obtain the column ID of $\bA$ as $\bA \approx \bm{C} \bm{P}$. The coefficient computation methods will be explained in Section \ref{subsec:coeff_compute}.
The detailed implementation of the online CSS is described in Algorithm \ref{alg:CSS_residual}. In Step 10, we notice that the selected columns $\bm{C}$ keep expanding without replacement and thus the algorithm may not be suitable for performing CSS with a fixed-size subset. 

\begin{algorithm}[H]
        \SetAlgoLined
        \KwIn{Data matrix $\bA_{m\times n}$, target rank $k$, parameter $\xi > 0 $}
        Initialize $\bm{C}_{\pre} = \emptyset$, $\bm{C}_{\current} = \emptyset$ and $\sigma = 0$ \\
        \While{A is not completely read through}{
        Read next column of $\bA$, denoted by $\bm{a}$, into RAM \\
        Let $p_a := \dfrac{k\Vert\Pi^{\perp}_{\bm{C}_{\pre}}\bm{a}\Vert_2^2}{\xi}$.\\
        With probability $\min(p_a,1)$, add $\bm{a}$ to $\bm{C}_{\current}$ \\

\If{$p_a<1$}{
    Increment $\sigma \leftarrow \sigma + p_a$ 
    }
\If{$p_a\ge1$ or $\sigma\ge 1$}{
    Set $\bm{C}_{\pre} \gets \bm{C}_{\pre} \cup \bm{C}_{\current}$ and reset $\sigma = 0$ and $\bm{C}_{\current} = \emptyset$ 
  }
        
        }
        \KwOut{$\bm{C} = \bm{C}_{\pre} \cup \bm{C}_{\current}$}
        \caption{Online CSS using residual-based sampling \cite{bhaskara2019residual}}
        \label{alg:CSS_residual}
\end{algorithm}

\subsubsection{Streaming ridge leverage score based CSS}
\label{subsubsec:ridge_leverage_CSS}
In \cite{cohen2017input}, Cohen \textit{et al.} proposed a streaming CSS algorithm via ridge leverage score sampling. The idea is to select columns using non-uniform probabilities determined by ridge leverage score. The leverage score of the $j$th column $\bm{a}_j$ of $\bA$ is defined as:
\begin{equation}
        \tau_j^0 = \bm{a}_j^T(\bA\bA^T)^{\pinv}\bm{a}_j,
        \label{eq:leverage-score}
\end{equation}
where $\tau_j^0$ measures how important $\bm{a}_j$ is in composing the range of $\bA$. It is maximized at 1 when $\bm{a}_j$ is linearly independent of the other columns of $\bA$ and decreases when many other columns approximately align with $\bm{a}_j$. In the low-rank approximation problem, the rank-$k$ leverage score is computed with respect to the rank-$k$ approximation $\bA_k$ of $\bA$ to evaluate how important each column is in composing the top $k$ singular directions of $\bA$'s range. However, $\bA_k$ is not always unique and the low-rank leverage score is sensitive to matrix perturbation which largely limits its ability to obtain accurate scores for sampling. To resolve this issue, the \emph{ridge} leverage score is defined by adding a regularization term to Eq.~(\ref{eq:leverage-score}),
\begin{equation}
        \tau_j = \bm{a}_j^T(\bA\bA^T + \lambda \bm{I})^{\pinv}\bm{a}_j,
        \label{eq:ridge-leverage-score}
\end{equation}
where we set $ \lambda = \|\bA-\bA_k\|^2_F/k$. Adding $\lambda\bm{I}$ to $\bA\bA^T$ alleviates the effect of smaller singular directions, leading to smaller sampling probabilities for them. As shown in \cite{cohen2017input}, given $\epsilon<1$ and with probability $1-\delta$, the ridge leverage score-based CSS satisfies the error bound
\begin{equation}
        \|\bA- (\bm{C}\bm{C}^{\pinv}\bA)_k\|^2_F \leq (1+\epsilon) \|\bA-\bA_k\|^2_F,
\end{equation}
where $\bm{C}$ is constructed by selecting $t = \mathcal{O}(k \log k +k \log(1/\delta)/\epsilon)$ columns from $\bA$.

In the streaming CSS settings, the storage limit for the loaded columns is assumed to be $t$. The basic idea of the algorithm is that whenever $t$ columns are loaded into RAM, the ridge leverage score for each column is computed and used to update the selected columns based on the score.

A convenient property of ridge leverage scores is that they are \emph{monotone} with respect to adding new columns, in the following sense:
\begin{lemma}
        \emph{\cite{cohen2017input}} (Monotonicity of ridge leverage score). For any $\bA \in \R^{m \times n}$ and vector $\bm{b} \in \R^{m}$, for every $j \in 1,2,...,n$ we have
        \begin{equation}
                \tau_j(\bA) \leq \tau_j(\bA\cup \bm{b}),
        \end{equation}
        where $\bA\cup \bm{b}$ denotes the column $\bm{b}$ appended to $\bA$ as the final column.
        \label{lemma:monotonicity}
\end{lemma}
Based on the monotonicity of the ridge leverage score (Lemma \ref{lemma:monotonicity}), whenever a new column is added to the original matrix, the ridge leverage scores of the already observed columns only decrease, which ensures that the previously selected columns are not biased with larger scores due to their earlier observation and thus newly observed columns can also be selected into the column subset. To ensure this monotonicity property holds for the low-rank leverage scores, where $\lambda$ depends on $\bA$, a generalized ridge leverage score of $\bA \in \R^{m \times n}$ with respect to any matrix $\bm{B} \in \R^{m \times n_B}$ is defined in \cite{cohen2017input} as 
\begin{equation}
        \tau_i^{\bm{B}}(\bA) = 
        \begin{cases}

                        \bm{a}_i^T\left(\bm{B}\bm{B}^T+ \dfrac{\|\bm{B}-\bm{B}_k\|^2_F}{k}\bm{I}\right)^{\pinv}\bm{a}_i &  \text{for } \bm{a}_i \in \spantext\left(\bm{B}\bm{B}^T+ \dfrac{\|\bm{B}-\bm{B}_k\|^2_F}{k}\bm{I}\right), \\
                        \infty   & \text{otherwise}.
        \end{cases}
        \label{eq:g_ridge_leverage_score}
\end{equation}
The generalized monotonicity bound of Eq.~(\ref{eq:g_ridge_leverage_score}) has been proven in \cite{cohen2017input}. 
In the streaming CSS process, $\bm{B}$ is chosen as the {\it frequent-directions} \cite{ghashami2016frequent} sketch of $\bm{A}$. The frequent-directions algorithm maintains the invariant that the last column of $\bm{B}$ is always all-zero valued. When a new column from the input matrix $\bA$ is observed, it replaces the all-zero valued column of $\bm{B}$. Then, the last column of $\bm{B}$ is nullified by a two-stage process. First, the sketch is rotated using its SVD such that its columns are orthogonal and in descending magnitude order. Then, the sketch columns norms are “shrunk” so that at least one of them is set to zero. The algorithm  produces a $\bm{B}$ satisfying
\begin{equation}
    \bm{BB}^T \prec \bm{AA}^T \quad \texttt{and}\qquad \|\bm{AA}^T - \bm{BB}^T\|_2 \leq 2\|\bA\|^2_F/n_B.
\end{equation}

At the beginning of the CSS process, $\bm{C} = \bm{0}_{m\times t}$, $\bm{D} = \bm{0}_{m\times t}$ matrices are initialized to store currently selected columns and new data columns, respectively. The vector $[\tau_1^{old},...,\tau_t^{old} ] = \bm{1}$ is also initialized to store the ridge leverage score of $\bm{C}$ with respect to the already observed data columns. With the arrival of the data sequence, $\bm{D}$ and the frequent-directions sketch $\bm{B}$ are first updated using every new data column until $t$ columns are read. The ridge leverage scores 
$\bm{\tau}^{\bm{D}}$ for new columns in $\bm{D}$ as well as the ridge leverage scores $\bm{\tau}$ for the previously selected columns $\bm{C}$ are recalculated with respect to the updated sketch $\bm{B}$, cf.\ \texttt{ApproximateRidgeScores} in Algorithm~\ref{alg:ridge_leverage_score_CSS}. 
Then, the normalized leverage score is used as the probability for selecting new columns and updating $\bm{C}$. For the $i^\text{th}$ column in $\bm{C}$, if the ridge leverage score $\tau_i$ decreases compared to $\tau_i^{old}$, we will discard that column with a probability equal to the proportion that its ridge leverage score decreased by, otherwise, the column is kept in the subset. If the column is discarded, we will select a new column from $\bm{D}$ based on the ridge leverage score $\tau^{\bm{D}}$. This process is similar to reservoir sampling \cite{vitter1985random} in that at any point in the stream, we have a set of columns sampled using true ridge scores concerning the currently observed matrix as the universal criteria. The error and complexity of the streaming ridge leverage score-based CSS is proved in \cite{cohen2017input} as the following theorem:
\begin{theorem}
        (Streaming column subset selection) \emph{\cite{cohen2017input}} Algorithm~\ref{alg:ridge_leverage_score_CSS} achieves computing a $(1 + \epsilon)$ error column subset $\bm{C}$ with $t = \mathcal{O}(k \log k +k \log(1/\delta)/\epsilon)$ columns in a single-pass over $\bA$. The algorithm uses $\mathcal{O}(mk)$ space in addition to the space required to store the subset and succeeds with probability $1-\delta$.
\end{theorem}

\noindent\scalebox{.9}{
\begin{minipage}{1.03\textwidth}
\begin{algorithm}[H]
        \SetAlgoLined
        \KwIn{Data matrix $\bA_{m\times n}$, target rank $k$, column buffer size $t$, frequent-directions sketch column size $n_B = 3k$, accuracy $\epsilon$, success probability ($1-\delta$)}
        \KwOut{Selected columns $\bm{C}_{m\times t} $}
        Initialize the Frobenius norm of the observed data matrix as $\|\bA_{obs}\|^2_{F} = 0$\\
        Set $count = 0$, $\bm{B}= \bm{0}_{m\times n_B}$, $\bm{C} = \bm{0}_{m\times t}$, $\bm{D} = \bm{0}_{m\times t}$\\
        Set the index set for the basis columns and new columns as $\mathcal{J}^{\bm{C}}, \mathcal{J}^{\bm{D}} \gets \bm{0}_{t\times 1}$ \\
        Initialize $[\tau_1^{old},...,\tau_t^{old} ] = 1$ for storing ridge leverage scores of $\bm{C}$ with respect to previous read data columns.\\
        \For{$j=1,...,n$}{
        Read the next column $j$ of $\bA$ into RAM, denoted by $\bm{a}_j$ \\
        $\bm{B} = \texttt{FrequentDirections}(\bm{B}, \bm{a}_j)$ \tcp*{Or use FastFrequentDirections (cf. Algorithm 2 in \cite{ghashami2016frequent})}
        \eIf{$count < t$}{
        $\bm{d}_{count} = \bm{a}_j$, $\mathcal{J}^{\bm{D}}_{count} = j$, $\|\bA_{obs}\|^2_{F} \gets \|\bA_{obs}\|^2_{F} + \|\bm{a}_j\|_2^2$, $count \gets count + 1$ \\
        }
        {
        Update leverage scores of columns in $\bm{C}$ as $[\tau_1,...,\tau_t] = \min \{[\tau_1^{old},...,\tau_t^{old} ], \texttt{ApproximateRidgeScores}(\bm{B}, \bm{C}, \|\bA_{obs}\|^2_{F})\}$\\
        Approximate leverage scores of new columns in $\bm{D}$ as $[\tau_1^{\bm{D}},...,\tau_t^{\bm{D}} ] = \texttt{ApproximateRidgeScores}(\bm{B}, \bm{D},\|\bA_{obs}\|^2_{F})$\\
        \For{$i=1,...,t$}{
                \If{$\bm{c}_i \neq 0$}{
                        With probability $(1-\tau_i/\tau_i^{old})$, set $c_i = 0$, $\tau_i^{old} = 1$ and $\mathcal{J}^{\bm{C}}_i = 0$ \\
                        Otherwise set $\tau_i^{old} = \tau_i$\\
                }
                \If{$\bm{c}_i = 0$}{
                        \For{$r=1,...,t$}{
                                With probability $\dfrac{\tau_r^{\bm{D}} c (k \log k+k \log(1/\delta)/\epsilon)}{t}$, set $\bm{c}_i = \bm{d}_r$, $\tau_i^{old} = \tau_r^D$ and $\mathcal{J}^{\bm{C}}_i = \mathcal{J}^{\bm{D}}_r$\\
                        }
                }
        }
        $count = 0 $ \\
        }
        }
        \SetKwFunction{FMain}{FrequentDirections}
        \SetKwProg{Fn}{Function}{:}{}
        \Fn{\FMain{$\bm{B}$, $\bm{a}_j$}}{
        \tcp{Follow the idea in \cite{ghashami2016frequent} to update $\bm{B}$}
        Replace the $n_B^{\text{th}}$ all-zero valued column of $\bm{B}$: $\bm{B}(:,n_B) \gets \bm{a}_j$ \\
        Perform SVD on $\bm{B}$: $[\bm{U}, \bm{\Sigma}, \bm{V}] = \texttt{svd}(\bm{B})$ \tcp*{Thin SVD}
        $\bm{B} \gets \bm{U} \sqrt{\bm{\Sigma}^2 - \sigma_{n_B}^2 \bm{I}_{n_B}}$ \tcp*{$\sigma_{n_B}$ is the $n_B^{\text{th}}$ singular value of $\bm{B}$}
        \textbf{return} $\bm{B}$
        }
        \SetKwFunction{FMain}{ApproximateRidgeScores}
        \SetKwProg{Fn}{Function}{:}{}
        \Fn{\FMain{$\bm{B}$, $\bm{M}_{m\times t}$, $\|\bA_{obs}\|^2_{F}$}}{
        \tcp{$\bm{M}$ contains any $t$ columns from $\bA$.}
        \For{$i=1,...,t$}
        {
        $\tau_i = \bm{m}_i^T \left(\bm{BB}^T+\dfrac{\|\bA_{obs}\|^2_{F}-\|\bm{B}_k\|^2_F}{k}\bm{I}\right)^{\pinv} \bm{m}_i$
        }
        \textbf{return} $[\tau_1,...,\tau_t]$
        }
        \caption{Streaming CSS based on ridge leverage score \cite{cohen2017input}}
        \label{alg:ridge_leverage_score_CSS}
\end{algorithm}
\end{minipage}
}

The detailed implementation of the streaming ridge leverage score-based CSS is described in Algorithm \ref{alg:ridge_leverage_score_CSS}. Regarding the input target rank $k$ and column buffer size $t$, we set $k = t$ for simplicity in our implementation. Compared to the residual-based CSS explained in Section \ref{subsubsec:residual_CSS}, this method prunes and updates the selected columns $\bm{C}$ (Steps 12--27) and thus maintains a fixed RAM usage during the CSS process. Therefore, we developed our randomized ID algorithm based on the streaming ridge leverage score-based CSS method.
\FloatBarrier
\section{Online randomized ID based on streaming CSS algorithm}
\label{sec:Method}
\begin{figure}[bht]
        \centering
        \includegraphics[width = 0.8\textwidth]{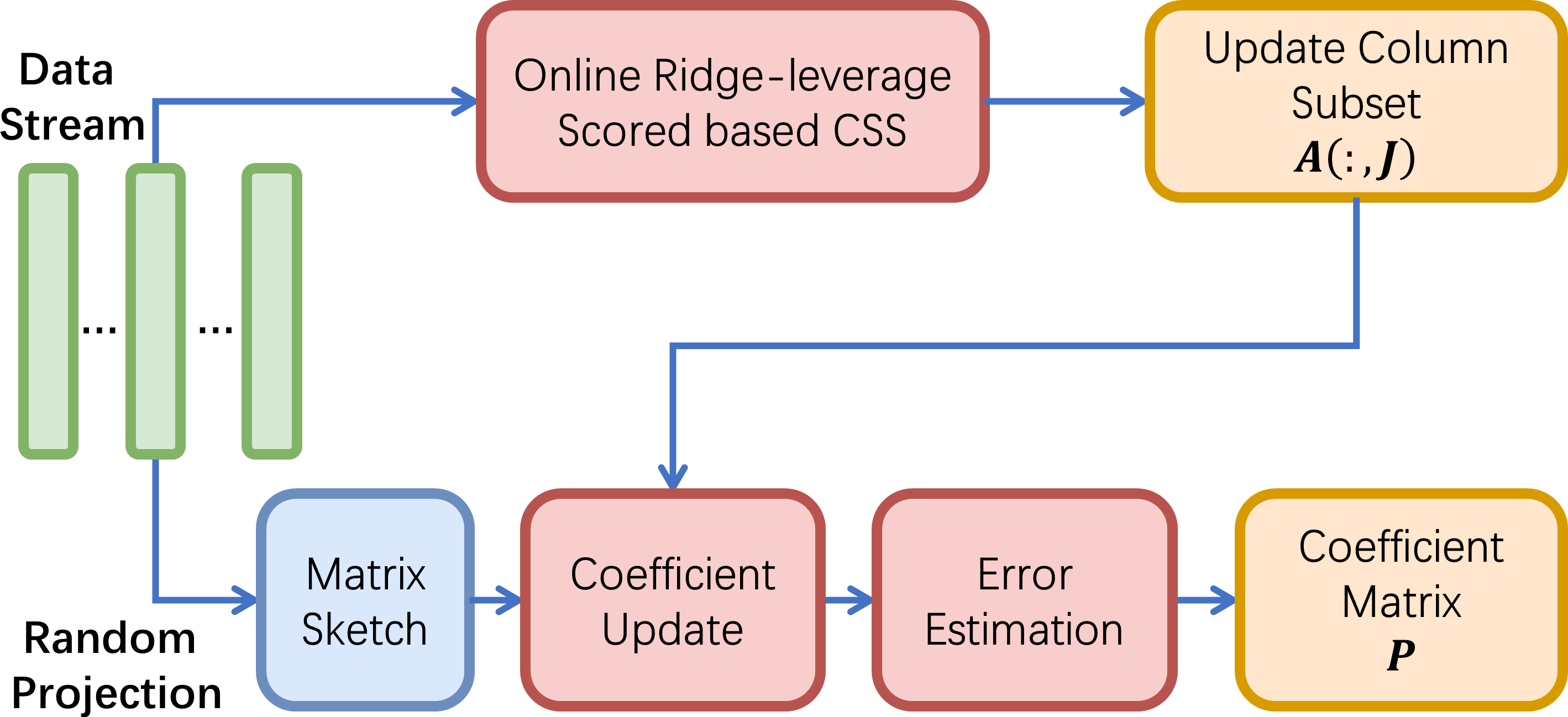}
        \caption{The workflow of our online randomized ID method.}
        \label{fig:workflow}
\end{figure}

The overall workflow of our online randomized ID approach is summarized in Fig.~\ref{fig:workflow} and explained below. There are two major procedures in the framework: (1) perform a single-pass column subset selection (CSS) to identify the appropriate basis columns $\bA_{\mathcal{J}}$ for approximation (Section \ref{subsec:ridge_leverage_CSS}), and (2) update the coefficient matrix $\bm{P}$ whenever the basis columns change (Section \ref{subsec:coeff_compute}). 
For the CSS procedure, we implement the streaming ridge leverage score-based CSS algorithm \cite{cohen2017input} to handle the streaming data. When a new data column arrives, the ridge leverage score of the new column is approximated using the currently selected columns, and the column subset is updated by random sampling based on the new leverage score. 

The standard coefficient update procedure calculates the coefficients by solving a least squares problem:
\begin{align}
        \bm{P} & = \underset{\bm{X}}{\argmin} \| \bA_{\mathcal{J}}\bm{X}- \bA\|_F^2 \label{eq:coeff_ID1}               \\
                   & = (\bA_{\mathcal{J}}^T\bA_{\mathcal{J}})^{-1}
                   \bA_{\mathcal{J}}^T \bA = 
                   \bA_{\mathcal{J}}^{\pinv} \bA.
        \label{eq:coeff_ID2}
\end{align}
However, the calculation of this $\bm{P}$ requires full access to the data matrix $\bA$, and hence requires a second pass over the data which we are trying to avoid. It is not possible to calculate the coefficients column-by-column during the first pass through the data since if a column $\bm{a}_i$ is loaded, the index set $\mathcal{J}$ may (eventually) contain indices $j$ that are greater than $i$ and so are not available at the time.

To address this issue, we apply random projections (cf.~\cite{halko2011finding}) to compute a small matrix sketch of $\bA$ and store it in memory for the coefficient computation; this can be done in the same pass as the online CSS algorithm. Since the matrix sketch is used to compute the coefficient and an optimal coefficient is not guaranteed, we use four different methods (Algs. \ref{alg:coeff_update1}--\ref{alg:coeff_update4}) to compute $\bm{P}$ simultaneously. To determine the best coefficients, we implement an error estimator using the NA-Hutch++ algorithm \cite{meyer2021hutch++, jiang2021optimal} to approximate the Frobenius norm errors of the reconstructed data matrices (Section \ref{subsec:hutch++}). The NA-Hutch++ is a single-pass stochastic algorithm that can reuse the matrix sketch obtained from random projection of $\bA$ for error estimation and thus can be easily accommodated in the online workflow. 
In addition, we apply a least squares method to estimate the gradient field, which is further used in the computation of the parameter matrix to ensure the gradient information of the physical field is well reconstructed (Section \ref{subsec:gradient_estimate}). 

\subsection{Streaming ridge leverage score based CSS using random projection}
\label{subsec:ridge_leverage_CSS}
In Section \ref{subsubsec:ridge_leverage_CSS}, the ridge leverage score-based streaming CSS algorithm utilizes the frequent-directions sketch for approximating the ridge leverage score (Eq.~\eqref{eq:g_ridge_leverage_score}). Since the data matrix $\bA \in \R^{m \times n}$ sometimes has very large dimension $m$, it can be impractical to compute $\bm{B}\bm{B}^T\in\R^{m\times m}$ due to its large dimension. Instead, we perform random projection on each column to collect a matrix sketch $\bS$ which can also be used to compute the coefficient matrix and estimate the reconstruction in our workflow. Here, we define the ridge leverage score with respect to $\bS$ as 
\begin{equation}
        \tau_i = (\bOmega \bm{a}_i)^T\left(\bS\bS^T+\dfrac{\|\bA_{obs}\|_F^2-\|\bm{S}_k\|^2_F}{k}\bm{I}\right)^{\pinv}(\bOmega \bm{a}_i),
        \label{eq:sketch_ridge_leverage_score}
\end{equation}
making it analogous to Eq.~\eqref{eq:g_ridge_leverage_score}. $\bOmega$ is the Gaussian random matrix; $\bA_{obs} = \bA(:,1:j)$ and $\bS=\bOmega\bA_{obs}$ are the observed $j$ columns of the data matrix and their sketch in the streaming setting, respectively; and $\bS_k$ is the best rank-$k$ approximation of $\bS$. To adopt the random projection in the CSS process, we make the following modifications based on the original streaming CSS algorithm (Algorithm \ref{alg:ridge_leverage_score_CSS}):
\begin{itemize}
        \item At the beginning of the CSS process, we generate the random projection matrix $\bOmega_{\ell \times m}$ and then initialize the sketch matrix $\bS = \bm{0}_{\ell \times n}$.
        \item The random projection is used to compute and collect the matrix sketch instead of frequent-directions.
        \item The ridge leverage scores for new columns in $\bm{D}$, $\bm{\tau}^{\bm{D}}$, and previously selected columns $\bm{C}$, $\bm{\tau}$, with respect to the updated sketch $\bS$ are computed using Eq.~\eqref{eq:sketch_ridge_leverage_score}. In particular, we substitute $\bm{a}_i$ in Eq.~\eqref{eq:sketch_ridge_leverage_score} with the column in $\bm{C}$ and $\bm{D}$ to compute the ridge leverage score of each column. 
\end{itemize}
The resulting pseudocode for the online randomized ID algorithm appears as Algorithm \ref{alg:rID_overview}, where the detailed implementation of the modified CSS method is described in Steps 9--27. In Steps 28--32, we compute and update the coefficient matrix for the ID with the help of an error estimator. The detailed algorithms of Steps 29 and 30 will be explained in Sections \ref{subsec:coeff_compute} and \ref{subsec:hutch++}. 

\noindent\scalebox{.9}{
\begin{minipage}{1.03\textwidth}
\begin{algorithm}[H]
        \SetAlgoLined
        \KwIn{Data matrix $\bA_{m\times n}$, target rank $k$, basis size $t$, accuracy $\epsilon$, success probability ($1-\delta$)}
        \KwOut{$\bA \approx \bA_\mathcal{J} \bm{P}$, where $\mathcal{J}$ contains the indices of the selected columns}
        Choose $\ell = k + p$, where $p$ is the oversampling parameter. \\
        Generate Gaussian random matrix $\bOmega = \texttt{randn}(\ell, m)$ \\
        Initialize the Frobenius norm of the observed data matrix as $\|\bA_{obs}\|^2_{F} = 0$\\
        Set $count = 0$, $\bm{C} = \bm{0}_{m\times t}$, $\bm{D} = \bm{0}_{m\times t}$ \\
        Set $\bS = \bm{0}_{\ell\times t}$, $\bm{SS}^T = \bm{0}_{\ell\times \ell}$\\
        Set the index set for the basis columns, the basis columns before the update, and the new columns as $\mathcal{J}^{\bm{C}}, \mathcal{J}^{\bm{C}}_{prev},\mathcal{J}^{\bm{D}} \gets \bm{0}_{t\times 1}$ \\
        Initialize $[\tau_1^{old},...,\tau_t^{old} ] = 1$ for storing ridge leverage scores of $\bm{C}$ with respect to previous read data columns.\\
        \For{$j=1,...,n$}{
        Read the next column $j$ of $\bA$ into RAM, denoted by $\bm{a}_j$ \tcp*{For analysis purposes, write $\bA_{obs}=\bA(:,1:j)$ to denote observed columns}
        Compute and store sketched data $\bm{s}_j = \bOmega \bm{a}_j$, $\bS \gets [\bS, \bm{s}_j]$ \\
        $\bm{SS}^T \gets \bm{SS}^T + \bm{s}_j \bm{s}_j^T$ \\
        \eIf{$count < t$}{
        $\bm{d}_{count} = \bm{a}_j$, $\mathcal{J}^{\bm{D}}_{count} = j$, $\|\bA_{obs}\|^2_{F} = \|\bA_{obs}\|^2_{F} + \|\bm{a}_j\|_2^2$, $count \gets count + 1$ \\
        }
        {
        \tcp{Follow the idea in \cite{cohen2017input} to compute ridge leverage score and prune column basis.}
        Update leverage scores of columns in $\bm{C}$ as $[\tau_1,...,\tau_t] = \min \{[\tau_1^{old},...,\tau_t^{old} ], \texttt{ApproximateRidgeScores\_Sketch}(\bm{SS}^T, \bOmega, \bm{C}, \|\bA_{obs}\|^2_{F})\}$\\
        Approximate leverage scores of new columns in $\bm{D}$ as $[\tau_1^{\bm{D}},...,\tau_t^{\bm{D}} ] = \texttt{ApproximateRidgeScores\_Sketch}(\bm{SS}^T, \bOmega, \bm{D},\|\bA_{obs}\|^2_{F})$\\
        \For{$i=1,...,t$}{
                \If{$\bm{c}_i \neq 0$}{
                        With probability $(1-\tau_i/\tau_i^{old})$, set $c_i = 0$, $\tau_i^{old} = 1$ and $\mathcal{J}^{\bm{C}}_i = 0$ \\
                        Otherwise set $\tau_i^{old} = \tau_i$\\
                }
                \If{$\bm{c}_i = 0$}{
                        \For{$r=1,...,t$}{
                                With probability $\dfrac{\tau_r^{\bm{D}} c (k \log k+k \log(1/\delta)/\epsilon)}{t}$, set $\bm{c}_i = \bm{d}_r$, $\tau_i^{old} = \tau_r^D$ and $\mathcal{J}^{\bm{C}}_i = \mathcal{J}^{\bm{D}}_r$\\
                        }
                }
        }
        
  \For{$q=4,\ldots,7$}
    {
        $\bm{P}_q = \texttt{CoefficientUpdate}(\mathcal{J}^{\bm{C}}, \bOmega\bA_{obs}, \mathcal{J}^{\bm{C}}_{prev}, \bm{P}_{prev}, \text{algo}=q)$ \tcp*{
        Via Algorithms \ref{alg:coeff_update1}--\ref{alg:coeff_update4}, $\mathcal{J}^{\bm{C}}_{prev}$ and $\bm{P}_{prev}$ are only required for $q=6,7$} 
        $e_q$ = approximation error estimated via single-pass Hutch++ (Algorithm~\ref{alg:NAhutchPP})
    }
      $\bm{P} = \bm{P}_{q^\star}$ where $q^\star = \argmin_{q\in\{4,\ldots,7\}} e_q$\\      
        $\mathcal{J}^{\bm{C}}_{prev} = \mathcal{J}^{\bm{C}}, count = 0 $ 
        }
        }
        \SetKwFunction{FMain}{ApproximateRidgeScores\_Sketch}
        \SetKwProg{Fn}{Function}{:}{}
        \Fn{\FMain{$\bm{SS}^T$, $\bOmega$, $\bm{M}_{m\times t}$, $\|\bA_{obs}\|^2_{F}$}}{
        \tcp{$\bm{M}$ contains any $t$ columns from $\bA$.}
        \For{i=1,...,t}
        {
        $\tau_i = (\bOmega \bm{m}_i)^T \left(\bm{SS}^T+\dfrac{\|\bA_{obs}\|^2_{F}-\|\bm{S}_k\|^2_F}{k}\bm{I}\right)^{\pinv} (\bOmega \bm{m}_i)$
        }
        \textbf{return} $[\tau_1,...,\tau_t]$
        }
        \caption{Single-pass randomized ID algorithm \emph{(novel algorithm)}}
        \label{alg:rID_overview}
\end{algorithm}
\end{minipage}
}

\subsection{Coefficient computation methods}
\label{subsec:coeff_compute}
Once we obtain the column basis via the streaming CSS algorithm, we can compute the coefficient matrix for the reconstruction of the original matrix. 
Due to the online setting, we need to compute an accurate coefficient matrix for the observed data using the information stored in RAM. Since we use random projection to construct the data sketch in RAM, a standard approach is to solve the sketched version of Eq.~\eqref{eq:coeff_ID1} and obtain the coefficients, which we will detail in Algorithm \ref{alg:coeff_update1}. 
However there are several other one-pass variants that can compute the coefficients, and 
we implement four different methods (Algorithms \ref{alg:coeff_update1}--\ref{alg:coeff_update4}) to simultaneously compute the coefficient matrices during the CSS process and then choose the best one afterward. 
We assume the column basis matrix $\bA_\mathcal{J} \in \R^{m \times t}$ contains $t$ columns, and the coefficient matrix $\bm{P} \in \R^{t \times n_{obs}}$ is computed based on the observed data matrix $\bA_{obs} = \bA(:,1:n_{obs}) \in \R^{m \times n_{obs}}$, where $n_{obs}$ is the number of the observed columns. 

\subsubsection{Algorithm \ref{alg:coeff_update1}}
In Algorithm \ref{alg:coeff_update1}, we first explain the standard approach to obtain the coefficients using the random projection sketch. Whenever the column basis is updated, we first get the new column basis sketch $\bOmega\bA_{\mathcal{J}}$ from the data matrix sketch $\bOmega\bA_{obs}$. Then, we solve the full sketching least squares problem to compute the updated coefficient matrix as
\begin{equation}
        \bm{P} = \underset{\bm{X}}{\argmin} \| \bOmega\bA_{\mathcal{J}}\bm{X}- \bOmega\bA_{obs}\|_F^2 = (\bOmega\bA_{\mathcal{J}})^{\pinv} (\bOmega\bA_{obs} ).
        \label{eq:coeff_rID_fullsketch}
\end{equation}
By itself, this computation could be done once at the end of the simulation $n_{obs} = n$ since it will only improve as $n_obs$ increases, but we perform the computation at every iteration since having a preliminary estimate of the coefficients will be useful for some of the other coefficient update methods.
\begin{algorithm}[H]
        \SetAlgoLined
        \KwIn{Index set $\mathcal{J}$ for the updated basis columns, current data matrix sketch $\bOmega \bA_{obs}$}
        Solve the least squares problem $\bm{P} = \underset{\bm{X}}{\argmin} \|\bOmega \bA_\mathcal{J}\bm{X}- \bOmega \bA_{obs} \|_F^2$ \\
        \KwOut{$\bm{P}_{obs}$}
        \caption{Coefficient computation by directly solving full sketching least squares problem}
        \label{alg:coeff_update1}
\end{algorithm}

\subsubsection{Algorithm \ref{alg:coeff_update2}}
To improve the accuracy of the full sketching method, we can take advantage of the fact we have the unsketched column basis matrix $\bA_{\mathcal{J}}$ containing the original data for the coefficient computation. In Algorithm \ref{alg:coeff_update2}, 
we use both original and sketched column bases,
following similar approaches in \cite{pilanci2016iterative,becker2017robust},
setting
\begin{equation}
\bm{P}=
\bA_\mathcal{J}^\pinv \bOmega^T\bOmega \bA = 
\argmin_{\bm{X}}
 \| \bA_\mathcal{J}\bm{X} \|_F^2
 - 2 \tr ( \bm{X}^T\bA_\mathcal{J}^T\bOmega^T (\bOmega\bA_{obs}) ) + \|\bA_{obs}\|_F^2 
\end{equation}
motivated by the expansion
\begin{equation}
 \| \bA_\mathcal{J}\bm{X}- \bA_{obs} \|_F^2
 = \| \bA_\mathcal{J}\bm{X} \|_F^2
 - 2 \tr ( \bm{X}^T\bA_\mathcal{J}^T\bA_{obs} ) + \|\bA_{obs}\|_F^2.
\end{equation}
or equivalently by noting that in Eq.~\eqref{eq:coeff_ID2} there is no need to replace $\bA_{\mathcal{J}}^{\pinv}$ with its sketched version.
As in Algo.~\ref{alg:coeff_update1}, this computation could be done once when $n_{obs}=n$ but we perform it every iteration due to the interaction with other methods.

\begin{algorithm}[H]
        \SetAlgoLined
        \KwIn{Index set $\mathcal{J}$ for the updated basis columns, current data matrix sketch $\bOmega \bA_{obs}$}
        Compute and store $\bm{G} = \bA_{\mathcal{J}}^T \bA_{\mathcal{J}}$, $\bm{Y} = (\bOmega\bA_{\mathcal{J}})^T (\bOmega \bA_{obs})$
        \tcp*{Approximate $\bm{Y} = \bA_{\mathcal{J}}^T \bA_{obs}$ using sketched data.}
        $\bm{P}_{obs} = \bm{G}^{-1} \bm{Y}$ \\
        \KwOut{$\bm{P}_{obs}$}
        \caption{Coefficient computation based on partial sketching of new data column}
        \label{alg:coeff_update2}
\end{algorithm}

\subsubsection{Algorithm \ref{alg:coeff_update3}}
\newcommand{\cJ}{\mathcal{J}}
At step $j$ of the loop, when column $\bm{a}_j$ is loaded, coefficients for $i\in\cJ$ with $i\le j$ can be directly calculated via Eq.~\eqref{eq:coeff_ID1}; the main difficulty is that at later steps, $\cJ$ will be increased to contain more indices, and $\bm{a}_j$ will no longer be available (hence the reason for sketching in Algorithms \ref{alg:coeff_update1} and \ref{alg:coeff_update2}).
In Algorithm \ref{alg:coeff_update3}, we exploit this by treating indices $i\le j$ differently than indices $i>j$.

Specifically, we re-use the old coefficients (corresponding to old basis elements that have not been pruned) to subtract off a partial estimate and leave behind a residual, and then update coefficients for new basis elements ($i>j$) via sketching (as in Algo.~\ref{alg:coeff_update1}) on the residual.
%
In detail, whenever the column basis is updated, we first determine the indices of the unchanged column basis $\mathcal{I}_{uc}$ in the previous index set $\mathcal{J}_{prev}$.
In the previous coefficient matrix, we zero out the rows corresponding to the replaced columns,
\begin{equation}
        \bm{P}_{prev}(i,:) \leftarrow
        \begin{cases}
                \bm{P}_{prev}(i, :), & i \in \mathcal{I}_{uc} \\
                \bm{0},              & \text{otherwise}
        \end{cases}
        \text{    for $i= 1,...,t$.}
        \label{eq:coeff_update_residual_zeroout}
\end{equation}
We then compute the approximation residual using the unchanged basis as
\begin{equation}
        \bm{r}_{A} = \bOmega \bA - \bOmega \bA_{\mathcal{J}} \bm{P}_{prev}.
\end{equation}
Then, we use the new column basis to approximate the residual and obtain the new coefficient matrix as
\begin{align}
        \delta\bm{P} & = \underset{\bm{X}}{\argmin} \|\bOmega\bA_\mathcal{J}\bm{X}- \bm{r}_A\|_F^2, \\
        \bm{P}       & = \bm{P}_{prev} + \delta\bm{P}.
\end{align}

\begin{algorithm}[H]
        \SetAlgoLined
        \KwIn{Index sets $\mathcal{J}$, $\mathcal{J}_{prev}$ for the updated and previous basis columns, respectively, the previous coefficient matrix $\bm{P}_{prev}$, current data matrix sketch $\bOmega \bA_{obs}$}
        Determine the unchanged column indices in $\mathcal{J}_{prev}$ as $\mathcal{I}_{uc}$,\\
        Update $\bm{P}_{prev}$ with $\mathcal{I}_{uc}$ using Eq.~\eqref{eq:coeff_update_residual_zeroout}. \\
        Compute approximation residual $\bm{r}_{A} = \bOmega\bA_{obs} - \bOmega\bA_{\mathcal{J}} \bm{P}_{prev}$ \\
        Solve the least squares problem $\delta\bm{P} =\underset{\bm{X}}{\argmin} \|\bOmega\bA_\mathcal{J}\bm{X}- \bm{r}_{A}\|_F^2$\\
        $\bm{P}_{obs} = \bm{P}_{prev} + \delta\bm{P}$\\
        \KwOut{$\bm{P}_{obs}$}
        \caption{Coefficient computation based on the residual of fully sketched least squares problem}
        \label{alg:coeff_update3}
\end{algorithm}

\subsubsection{Algorithm \ref{alg:coeff_update4}}
Since the column basis keeps updating with more representative columns during the online CSS process, we assume that the old column basis can always be represented using the new basis. Therefore, we propose to compute a transformation between the old and new basis for updating the old coefficients. In Algorithm \ref{alg:coeff_update4}, we utilize the QR decomposition of the selected column basis to compute the transformation. Whenever the column basis is updated, we compute the QR on the updated column basis
\begin{equation}
        \bA_{\mathcal{J}} = \bm{QR}.
\end{equation}
We then compute the transformation matrix from new to old column basis and update the old coefficient matrix as
\begin{equation}
        \bm{T} = \bm{R}^{\pinv}\bm{Q}^{T} \bA_{\mathcal{J}_{prev}}, \\
        \bm{P} = \bm{T} \bm{P}_{prev}.
\end{equation}
\begin{algorithm}[H]
        \SetAlgoLined
        \KwIn{Index sets $\mathcal{J}$, $\mathcal{J}_{prev}$ for the updated and previous basis columns, respectively, the previous coefficient matrix $\bm{P}_{prev}$, current data matrix sketch $\bOmega \bA_{obs}$}
        Compute $\bm{Q}, \bm{R} = qr(\bA_\mathcal{J}) $ \\
        Compute the transformation matrix from new to old column basis as $\bm{T} = \bm{R}^{\pinv}\bm{Q}^{T} \bA_{\mathcal{J}_{prev}}$ \\
        Update the coefficient matrix $\bm{P} = \bm{T} \bm{P}_{prev}$ \\

        \KwOut{$\bm{P}$}
        \caption{Coefficient computation based on QR decomposition update}
        \label{alg:coeff_update4}
\end{algorithm}

\subsection{Single-pass Hutch++ for error estimation}
\label{subsec:hutch++}
\newcommand{\omegaH}{\bOmega_1}
\newcommand{\omegaR}{\bOmega_2}
\newcommand{\omegaG}{\bOmega_3}
We seek an {\em a posteriori} estimate of the reconstruction error because it is inherently useful and because we can apply all 
four algorithms presented in Section~\ref{subsec:coeff_compute} simultaneously (using the same sketched data) and use the error estimate to select whichever method is best.
Given our emphasis on online compression, we adopt the NA-Hutch++ trace estimation algorithm in \cite{meyer2021hutch++} as it requires only a single pass over the data to estimate the Frobenius norm error $\|\bE\|_F = \|\bA - \bA_{\mathcal{J}}\bm{P}\|_F$. In more detail, 
the NA-Hutch++ algorithm builds on the well-known Hutchinson's stochastic trace estimator \cite{hutchinson1989stochastic}. Let $\bm{B} = \bE\bE^T \in \R^{m \times m}$ be a square matrix and $\bOmega=[\bm{\omega}_1, ...,\bm{\omega}_\ell]^T \in \R^{\ell \times m} $ a matrix containing i.i.d.\ random variables with mean 0 and variance 1. For each $\boldsymbol{\omega}_j \in \R^m$, the expectation of $\boldsymbol{\omega}_j^T\bm{B}\boldsymbol{\omega}_j$ equals the trace of $\bm{B}$. The Hutchinson estimator approximates $\|\bE\|_F^2=\tr(\bm{B})$ as
\begin{equation}
        \tr(\bm{B}) \approx \dfrac{1}{\ell}\sum_{j=1}^{\ell}\boldsymbol{\omega}_j^T\bm{B}\boldsymbol{\omega}_j = \dfrac{1}{\ell} \tr(\bOmega\bm{B}\bOmega^T),
        \label{eq:Hutchinson}
\end{equation}
where $\ell$ represents the number of ``matrix-vector'' multiplication queries. 
This amounts to a Monte Carlo approach, so accuracy improves very slowly as $\ell$ increases, which has motivated modern improvements, such as Hutch++ \cite{meyer2021hutch++}, which use low-rank approximation as control variates to improve the estimate, though these methods are not one-pass. Instead, we turn to the NA-Hutch++ algorithm which is one-pass and has nearly the same theoretical guarantees as Hutch++ (proved in \cite{meyer2021hutch++}), leading to more accurate trace estimation as compared to Eq.~\eqref{eq:Hutchinson}.
The NA-Hutch++ method splits the random sampling matrix $\bOmega$ into three  random matrices $\omegaH \in \R^{c_1\ell \times m}$, $\omegaR \in \R^{c_2\ell \times m}$ and $\omegaG \in \R^{c_3\ell \times m}$, where $c_1\ell,c_2\ell$, and $c_3\ell$ are all integers satisfying $(c_1+c_2+c_3) = 1$. Then, if the random matrices are centered iid sub-Gaussian,  given $c_2 > c_1$ and $\ell = \mathcal{O}(k \log(1/\delta))$, with probability $1-\delta$, the following approximation of $\bm{B}$
\begin{equation}
        \Tilde{\bm{B}} = (\omegaR\bm{B})^T(\omegaH\bm{B}\omegaR^T)^{\pinv}(\omegaH\bm{B}),
\end{equation}
satisfies $\|\bm{B} - \Tilde{\bm{B}}\|_F \leq 2 \|\bm{B} - {\bm{B}}_k\|_F$ \cite{clarkson2009numerical}, and thus is an accurate, albeit biased, estimate of $\tr(\bm{B})$,
and note that $\tr(\Tilde{\bm{B}})$ can be computed efficiently without explicitly constructing  $\Tilde{\bm{B}}$ by exploiting the cyclic properties of trace.
The final NA-Hutch++ 
estimate debiases the $\Tilde{\bm{B}}$ estimate by using it as a control variate within the Monte Carlo estimate, i.e.,
\begin{equation}
                \tr(\bm{B}) \approx 
                \underbrace{\tr((\omegaH\bm{B}\omegaR^T)^{\pinv}(\omegaH\bm{B})(\omegaR\bm{B})^T)}_{\tr(\Tilde{\bm{B}})} 
                 + \dfrac{1}{c_3\ell}[\tr(\omegaG\bm{B}\omegaG^T) - \tr(\omegaG\Tilde{\bm{B}}\omegaG^T)].
        \label{eq:NAhutchPP1}
\end{equation}
Note that we can sample the row of the random projection matrix $\bOmega$ to obtain random matrices $\omegaH$, $\omegaR$ and $\omegaG$ without duplicated rows. Accordingly, we can also obtain three query matrices $\omegaH\bE = \omegaH(\bA - \bA_{\mathcal{J}}\bm{P})$, $\omegaR\bE = \omegaR(\bA - \bA_{\mathcal{J}}\bm{P})$ and $\omegaG\bE = \omegaG(\bA - \bA_{\mathcal{J}}\bm{P})$ by reusing the corresponding rows of the sketched data $\bOmega(\bA - \bA_{\mathcal{J}}\bm{P})$, and thus the query process can be easily accommodated in the online randomized ID pipeline (Step 30 in Algorithm \ref{alg:rID_overview}).

However, we cannot simply substitute $\bm{B} = \bE\bE^T$ into  Eq.~\eqref{eq:NAhutchPP1} because that would require
access to 
the full matrix of $\bE$ when computing $\omegaH\bE\bE^T$ and $\omegaR\bE\bE^T$, which is impossible in the online pipeline. To address this issue, we adopt the expansion of $\|\bE\|_F^2$ as 
\begin{equation}
        \|\bE\|_F^2 = \|\bA - \bA_{\mathcal{J}}\bm{P}\|^2_F = \|\bA\|^2_F - 2 \tr(\bA(\bA_{\mathcal{J}}\bm{P})^T) + \underbrace{\tr((\bA_{\mathcal{J}}^T\bA_{\mathcal{J}})\bm{P}\bm{P}^T)}_{\|\bA_{\mathcal{J}}\bm{P}\|^2_F},
        \label{eq:NAhutchPP_expansion}
\end{equation}
where $\|\bA\|_F^2$ is always accessible in the streaming mode (Step 13 in Algorithm \ref{alg:rID_overview}); $\tr((\bA_{\mathcal{J}}^T\bA_{\mathcal{J}})\bm{P}\bm{P}^T)$ can be directly evaluated; and now $\tr(\bA(\bA_{\mathcal{J}}\bm{P})^T)$ can be estimated using Eq.~\eqref{eq:NAhutchPP1} as 
\begin{align}
\begin{split}
        \tr(\bA(\bA_{\mathcal{J}}\bm{P})^T) &\approx \tr((\omegaH\bA(\omegaR\bA_{\mathcal{J}}\bm{P})^T)^{\pinv}(\omegaH\bA\bm{P}^T)(\bA_{\mathcal{J}}^T\bA_{\mathcal{J}})(\omegaR\bA\bm{P}^T)^T) \\
        &+ \dfrac{1}{c_3\ell}\bigg[{\tr((\omegaG\bA)(\omegaG\bA_{\mathcal{J}}\bm{P})^T)} \\
        &- \tr((\omegaG\bA_{\mathcal{J}}\bm{P})(\omegaR\bA)^T(\omegaH\bA(\omegaR\bA_{\mathcal{J}}\bm{P})^T)^{\pinv}(\omegaH\bA)(\omegaG\bA_{\mathcal{J}}\bm{P})^T)\bigg].
\end{split}
\label{eq:NAhutchPP_2nd_term}
\end{align}
We can reuse the corresponding rows of the sketched data $\bOmega\bA$ and $\bOmega\bA_\mathcal{J}$ to obtain the query matrices $\omegaH\bA$, $\omegaR\bA$, $\omegaG\bA$, $\omegaH\bA_\mathcal{J}$, $\omegaR\bA_\mathcal{J}$ and $\omegaG\bA_\mathcal{J}$. $\bA_{\mathcal{J}}^T\bA_{\mathcal{J}}$ is evaluated and stored in Algorithm \ref{alg:coeff_update2} and thus can be reused for computing the last two terms in Eq.~\eqref{eq:NAhutchPP_expansion}. As a result, whenever the column basis $\bA_\mathcal{J}$ is updated, we only need the current data matrix sketch $\bOmega\bA_{obs}$ to evaluate the Frobenius norm error of different coefficient updating methods. The detailed implementation of the error estimator is shown in Algorithm \ref{alg:NAhutchPP}.

\newcommand{\indH}{\mathcal{I}_1}
\newcommand{\indR}{\mathcal{I}_2}
\newcommand{\indG}{\mathcal{I}_3}
\begin{algorithm}[H]
        \SetAlgoLined
        \KwIn{Current data matrix sketch $\bOmega\bA_{obs}$, 
        current Frobenius norm $\|\bA_{obs}\|_F^2$, selected basis columns $\bA_{\mathcal{J}}$, parameter matrix $\bm{P}$, random matrix $\bOmega_{\ell \times m}$, constants $c_1,c_2$, and $c_3$ where $c_1<c_2$ and $c_1+c_2+c_3 = 1$, and $c_1\ell, c_2\ell, c_3\ell \in \mathbb{N}$,}
        Compute and store $\bm{G} = \bA_{\mathcal{J}}^T\bA_{\mathcal{J}}$\\
        Compute $\|\bA_\mathcal{J}\bm{P}\|^2_F = \tr(\bm{G}\bm{P}\bm{P}^T)$ \\
        Sampling from the row of $\bOmega$ and obtain three index sets $\indH,\indR,\indG$ with sizes $c_1\ell, c_2\ell, c_3\ell$.\\
        Set sample matrices $\omegaH = \bOmega(\indH, :)$, $\omegaR = \bOmega(\indR, :)$, $\omegaG = \bOmega(\indG, :)$\\
        Obtain the query matrices $\omegaH\bA_{obs}$, $\omegaR\bA_{obs}$, $\omegaG\bA_{obs}$ from $\bOmega\bA_{obs}$, and $\omegaH\bA_\mathcal{J}$, $\omegaR\bA_\mathcal{J}$, $\omegaG\bA_\mathcal{J}$ from $\bOmega\bA_{\mathcal{J}}$ \\
        Estimate $\tr(\bA_{obs}(\bA_{\mathcal{J}}\bm{P})^T)$ using Eq.~\eqref{eq:NAhutchPP_2nd_term} \\
        Estimate $\|\bA_{obs} - \bA_{\mathcal{J}}\bm{P}\|^2_F = \|\bA_{obs}\|^2_F - 2\tr(\bA_{obs}(\bA_{\mathcal{J}}\bm{P})^T) + \|\bA_\mathcal{J}\bm{P}\|^2_F$  \\
        \KwOut{$\|\bA_{obs} - \bA_{\mathcal{J}}\bm{P}\|_F$}
        \caption{NA-Hutch++ \cite{meyer2021hutch++} for the error estimation of online randomized ID}
        \label{alg:NAhutchPP}
\end{algorithm}

\subsection{Enhancing the randomized ID with estimated gradient information}\label{subsec:gradient_estimate}
In many scientific applications, the quantities of interest (QoIs) are obtained by postprocessing the solution to the governing equations. For instance, in computational fluid dynamics, the vorticity field is computed based on the gradient of the velocity field. 
When only the solution data is compressed, however, the quality of the QoIs generated from the reconstructed data may be poor. To address this issue, we consider the scenario where the QoI is the gradient of the solution, and present a modification of the online randomized ID approach that incorporates the gradient information in the compression process. In Section \ref{sec:Result}, we show that for a fixed basis size, the modified ID approach leads to considerably more accurate solution gradients at the expense of a slight decrease in the accuracy of the reconstructed solution. 

\newcommand{\indi}{q}
\newcommand{\indj}{r}
To account for the gradient information in our compression framework, we first apply a least squares method to estimate the gradient of the reconstructed data assuming a graph structure for the underlying mesh; such gradient estimates are often referred to as \emph{simplex gradients}~\cite{Bortz1998}. Let $\mathcal{G} = (\mathcal{V}, \mathcal{E})$ be a finite graph, where $\mathcal{V}$ and $\mathcal{E}$ denote the vertices and edges of the graph, respectively; $\bm{x}(\mathcal{V})$ denotes the $d$-dimensional spatial coordinates of the vertices; $F(\mathcal{V})$ denotes the physical field of interest defined on this graph, i.e., the data stored in $\bA$. For each vertex $v_\indi \in \mathcal{V}$, we define its 1-ring neighborhood $\mathcal{N}_\indi = \{v_\indj\}_{\indj=1}^{n_\indi}$ such that each vertex $v_\indj$ is directly connected to $v_\indi$ and $n_\indi$ is the number of the neighboring vertices, and denote $\mathcal{N}_\indi(\indj)$ as the $\indj^{\text{th}}$ vertex in the 1-ring neighborhood. In the 1-ring neighborhood of each node $v_\indi$, we have
\begin{equation}
        F(v_\indj) - F(v_\indi) \approx \nabla F(v_\indi) \cdot \left(\bm{x}(v_\indj) - \bm{x}(v_\indi)\right),\quad \indj= 1,...,n_\indi.
        \label{eq:grad_approx_nodal}
\end{equation}
By collecting all the $\bm{x}(v_\indj) - \bm{x}(v_\indi)$ in a matrix $\bm{K}_\indi$ and $F(v_\indj) - F(v_\indi)$ in a vector $\bm{f}_\indi$, we have 
\begin{equation}
        \bm{K}_\indi  = \begin{bmatrix}
                \bm{x}(v_1) - \bm{x}(v_\indi)\\
                 \vdots \\
                 \bm{x}(v_\indj) - \bm{x}(v_\indi)\\
                 \vdots \\
                 \bm{x}(v_{n_\indi}) - \bm{x}(v_\indi)\\
                \end{bmatrix}, \; 
        \bm{f}_\indi  = \begin{bmatrix}
                        F(v_{1}) - F(v_\indi) \\
                        \vdots\\
                        F(v_{\indj}) - F(v_\indi) \\
                         \vdots \\
                        F(v_{n_\indi}) - F(v_\indi)
                        \end{bmatrix}.
\end{equation}
We rewrite Eq.~\eqref{eq:grad_approx_nodal} as the linear system
\begin{equation}
        \bm{K}_\indi \nabla F(v_\indi) = \bm{f}_\indi,
\end{equation}
which is usually overdetermined ($n_\indi > d$) for the commonly used meshes such as the 2D or 3D structured grids. The least squares estimate  $\widehat{\nabla F}(v_\indi)$ of $\nabla F(v_\indi)$ is thus
\begin{equation}
        \widehat{\nabla F}(v_\indi) = \bm{K}_\indi^{\pinv} \bm{f}_\indi.
\end{equation}
By computing and assembling the local matrix $\bm{K}_\indi^{\pinv}$ of each vertex on the entire graph, we can build the global gradient estimation operator $\bm{G}$ as
\begin{equation}
        \bm{G}  = \begin{bmatrix}
                        \bm{G}^1 \\
                        \vdots\\
                        \bm{G}^p \\
                         \vdots \\
                        \bm{G}^d
                        \end{bmatrix},
\end{equation}
where $\bm{G}^p$ is the gradient estimation matrix of the $p^{\text{th}}$ spatial dimension and is defined as
\begin{equation}
        \bm{G}^p(\indi,\indj) \leftarrow
        \begin{cases}
                \bm{K}_{\indi}^{\pinv}(p, k), & v_\indj = \mathcal{N}_\indi(k) \\
                -\sum_{k=1}^{n_\indi} \bm{K}_{\indi}^{\pinv}(p,k), & v_\indj = v_\indi \\
                \bm{0},              & \text{otherwise}
        \end{cases}
        \quad\text{for $\indi= 1,...,m$,\  $\indj= 1,...,m$.}
                        \label{eq:gradient_estimation_mat}
\end{equation}
so that $\bm{G}\bm{a}$ is a vector of size $dm$ containing the gradient estimates for all $m$ nodes.
Note that $\bm{G}$ can be stored as a sparse matrix since each row has only $n_i+1$ non-zero elements. It can be used in both CSS and coefficient computation procedures for gradient estimation and improving the reconstruction of the gradient field. During the CSS procedure, the gradient information helps to select a better column basis accounting for the gradient information. In more detail, in Step 10 of the CSS algorithm (Algorithm \ref{alg:rID_overview}), we obtain the estimated gradient $\bm{G}\bm{a}_j$ for each new data column $\bm{a}_j$, and concatenate the estimated gradient vector with $\bm{a}_j$ to approximate the ridge leverage scores for column selection. The gradient estimate is discarded after the ridge leverage score is computed to reduce memory usage. The ID coefficient matrix can then be obtained by the augmented optimization problem
\begin{equation}
        \underset{\bm{X}}{\argmin}\left( \| \bA_{\mathcal{J}}\bm{X}- \bA\|^2_F + \lambda \sum_{p=1}^d \| \bm{G}^p\bA_{\mathcal{J}}\bm{X}- \bm{G}^p \bA\|^2_F\right),
        \label{eq:opt_grad}
\end{equation}
where $\lambda$ is a regularization parameter. In our implementation, 
we make this suitable for streaming by approximating the least-squares problem of Eq.~\eqref{eq:opt_grad} using full sketching:
\begin{equation}
        \bm{P}(\lambda) =\underset{\bm{X}}{\argmin}\left( \| \bOmega(\bA_{\mathcal{J}}\bm{X}- \bA)\|^2_F + \lambda \sum_{p=1}^d \| \bOmega(\bm{G}^p\bA_{\mathcal{J}}\bm{X}- \bm{G}^p \bA)\|^2_F\right).
        \label{eq:opt_grad_fullsketch}
\end{equation}
To determine the best $\lambda$, we follow the generalized cross-validation (GCV) method \cite{golub1979generalized} to minimize the GCV function
\begin{equation}
        \text{GCV}(\lambda) = \frac{\|(\bm{I}-\bm{C}^{\natural}(\lambda))\bA\|_F^2}{\left(\tr(\bm{I}-\bm{C}^{\natural}(\lambda))
        \vphantom{\bm{C}^{\natural}_{\bOmega}}
        \right)^2}
\end{equation}
where $\bm{C}^{\natural}(\lambda)$ is defined as
\begin{equation}
        \bm{C}^{\natural}(\lambda) = \bA_{\mathcal{J}}\left(\bA_{\mathcal{J}}^T\bA_{\mathcal{J}}+\lambda \sum_{p=1}^d (\bm{G}^p\bA_{\mathcal{J}})^T\bm{G}^p\bA_{\mathcal{J}}\right)^{-1}\left(\bA_{\mathcal{J}} + \lambda \sum_{p=1}^d (\bm{G}^p)^T\bm{G}^p\bA_{\mathcal{J}}\right)^T.
\end{equation}
such that it satisfies
\begin{equation}
        \bA_{\mathcal{J}} \bm{P}(\lambda) = \bm{C}^{\natural}(\lambda)\bA.
\end{equation}
In particular, we use the random projection matrix $\bOmega$ and the sketched data $\bOmega\bA$ to approximate the GCV function in a single pass as
\begin{equation}
        \text{GCV}(\lambda) \approx \frac{\|(\bm{I}-\bm{C}^{\natural}_{\bOmega}(\lambda))\bOmega\bA\|_F^2}{\left(\tr(\bm{I}-\bm{C}^{\natural}_{\bOmega}(\lambda))\right)^2},
\end{equation}
where
\begin{equation}
        \bm{C}^{\natural}_{\bOmega}(\lambda) = \bOmega\bA_{\mathcal{J}}\left((\bOmega\bA_{\mathcal{J}})^T(\bOmega\bA_{\mathcal{J}})+\lambda \sum_{p=1}^d (\bOmega\bm{G}^p\bA_{\mathcal{J}})^T\bOmega\bm{G}^p\bA_{\mathcal{J}}\right)^{-1}\left(\bOmega\bA_{\mathcal{J}} + \lambda \sum_{p=1}^d (\bOmega \bm{G}^p\bOmega^T)^T \bOmega \bm{G}^p\bA_{\mathcal{J}}\right)^T.
        \label{eq:GCV_COmega}
\end{equation}
To determine the optimal $\lambda$ value, we perform the golden section search \cite{kiefer1953sequential} within a range of $10^{-3}$ to $10^{3}$ at the end of compression. For different $\lambda$, we directly compute the GCV function using the sketched data. Specifically, since $\bA_{\mathcal{J}}$ is fixed at the end of compression, the terms related to $\bA_{\mathcal{J}}$ in Eq.~\eqref{eq:GCV_COmega} can be computed once and reused for different $\lambda$, making it computationally efficient to evaluate the GCV function. 

\section{Numerical experiments}
\label{sec:Result}
The online randomized ID algorithms proposed in this work are tested on the following three datasets:
\begin{itemize}
        \item The JHU turbulence flow data from the direct numerical simulation of a turbulent channel flow in a domain of size $8\pi \times 2 \times 3\pi$, using $2048 \times 512 \times 1536$ nodes~\cite{perlman2007data,li2008public,graham2016web}. The dataset contains 4,000 frames of data after the simulation reaches a statistically stationary state. The frictional Reynolds number is $Re_\tau \approx 1000$.
        \item The ignition dataset consists of species concentration in a jet diffusion simulation computed on a $50 \times 50$ grid over 450 time steps.
        \item The NSTX Gas Puff Image (GPI) data \cite{zhao2020sdrbench, SDRBench_git} consists of 2D time-series data on $80 \times 64$ grid in $300,\!000$ time steps.
\end{itemize}
In our experiments, we measure the relative Frobenius error 
of our algorithm and benchmark against the offline SVD and two-pass ID. We will demonstrate the effectiveness of our approach by comparing it with other approaches in terms of both CSS and coefficient calculation steps.

These datasets are not so large as to actually require a single-pass implementation on modern hardware, which allows us to compare to more alternative methods, but this also means that runtime comparisons are not meaningful. 
We do not report in detail on the wall-clock runtime of our algorithm in comparison to benchmark algorithms because this raises issues of implementation details and furthermore 
our main focus is on the streaming scenario. 
However, we do note that in our Python implementation, our algorithm had roughly similar performance time to the randomized SVD. For example, the turbulent channel flow with $256\times 256$ and rank 100 took $10.3$ seconds with the streaming ID compared to $16.5$ seconds for the randomized SVD; for the ignition data, all methods took under 3 seconds; and for the $80 \times 64 \times 200,000$ rank 100 gas puff data, the rID took 96 seconds compared to 80 for the SVD.

\subsection{Turbulent channel flow}
\label{subsec:channelflow}
\begin{figure}[!htbp]
        \centering
        \includegraphics[width =\textwidth]{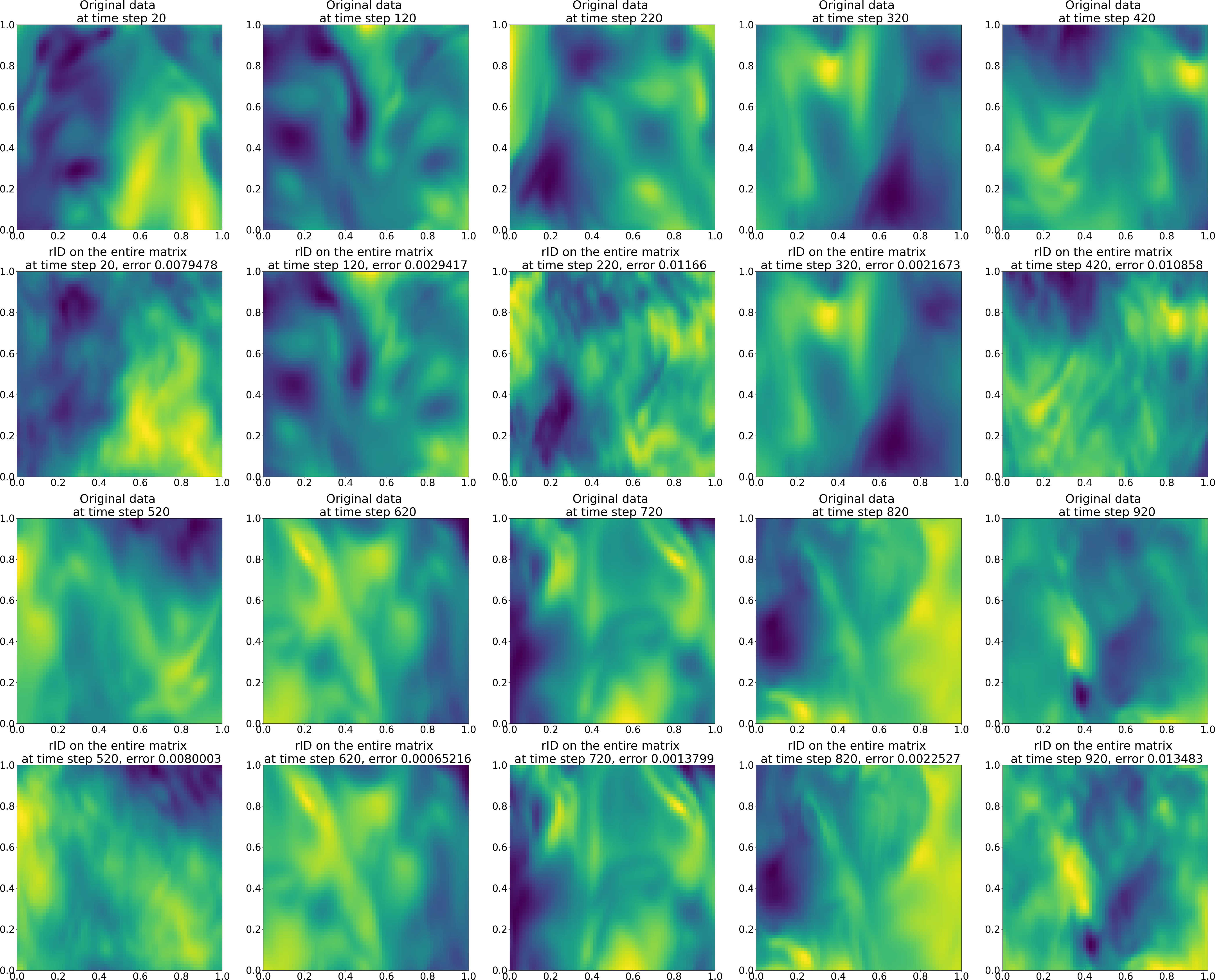}
        \caption{Rank $k = 50$ reconstruction of turbulence flow data over a $64 \times 64$ grid at different time steps. For each two rows, the top one represents the original data while the bottom one represents the reconstruction data.}
        \label{fig:channelflow_1}
\end{figure}

\begin{figure}[!htbp]
        \centering
        \includegraphics[width =\textwidth]{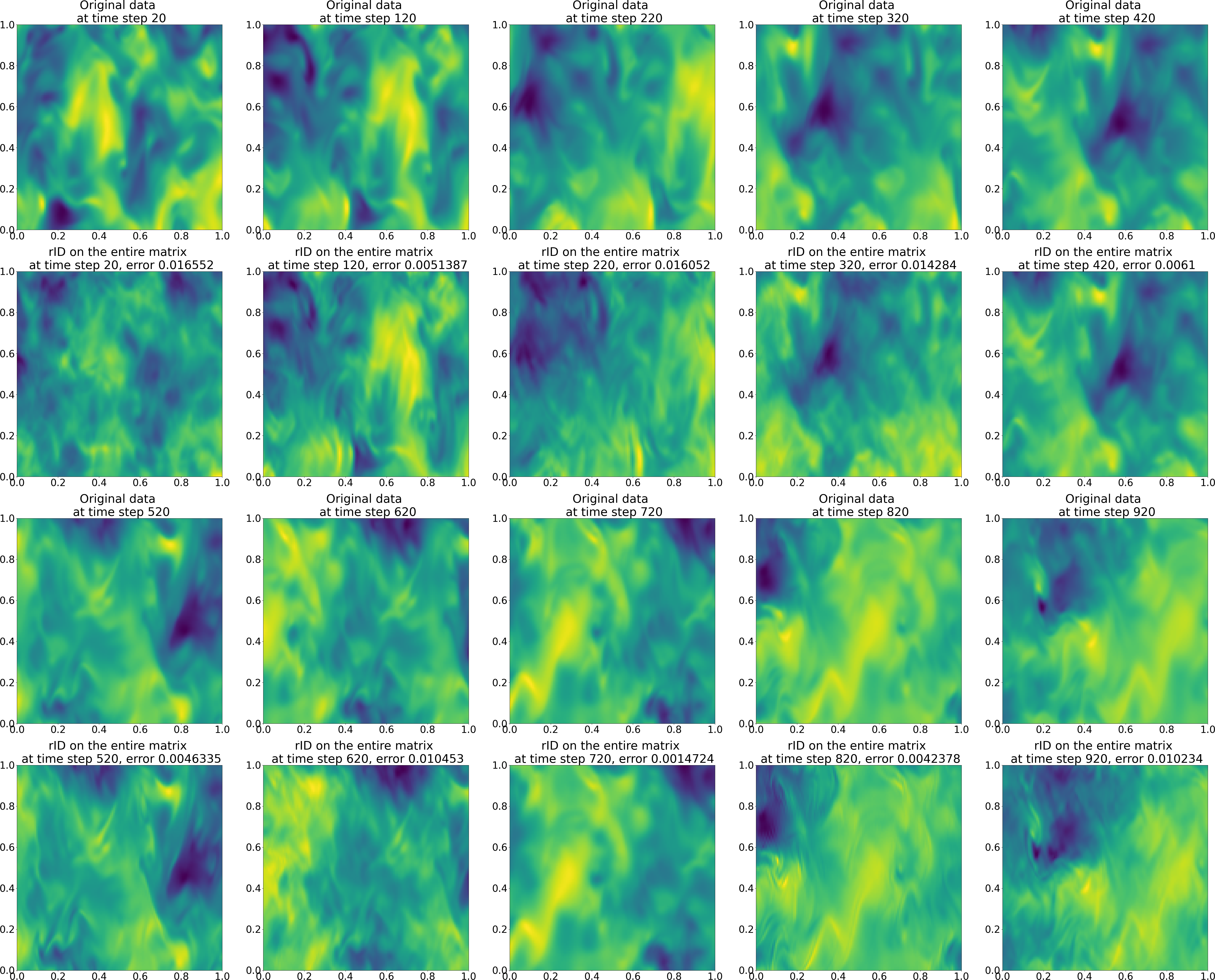}
        \caption{Rank $k = 50$ reconstruction of turbulence flow data over a $128 \times 128$ grid at different time steps. For each two rows, the top one represents the original data while the bottom one represents the reconstruction data.}
        \label{fig:channelflow_2}
\end{figure}

\begin{figure}[!htbp]
        \centering
        \includegraphics[width =\textwidth]{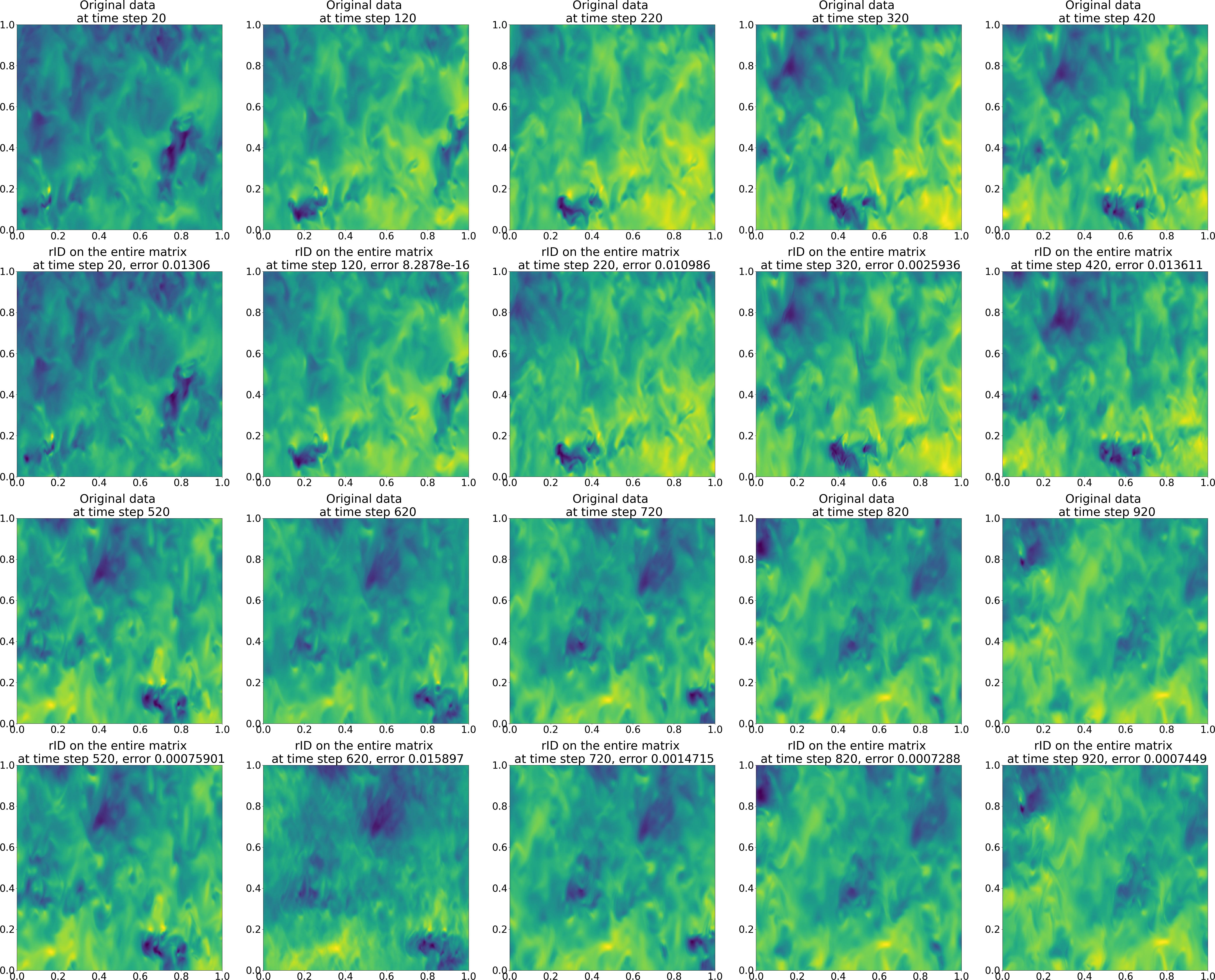}
        \caption{Rank $k = 50$ reconstruction of turbulence flow data over a $256 \times 256$ grid at different time steps. For each two rows, the top one represents the original data while the bottom one represents the reconstruction data.}
        \label{fig:channelflow_3}
\end{figure}

We first test our algorithm on the dataset extracted from the turbulence channel flow simulation result of 4000 frames \cite{perlman2007data,li2008public,graham2016web}. Regarding the boundary conditions of the channel flow, the periodic boundary conditions are imposed in the longitudinal and transverse directions, and the no-slip conditions are applied to the top and bottom walls of the channel.
We collect the $x$-direction velocity field ($u_x$) on $64 \times 64$, $128 \times 128$, and $256 \times 256$ $x$-$z$ grids at the $y=1$ plane through time, where the $x$-axis is defined as along the channel, the $y$-axis is perpendicular to the top and bottom walls of the channel, and the $x$-$z$ plane is parallel to the wall. 
Therefore, we obtain three data matrices of 4000 snapshots with 4096, 16,384, and 65,536 grid points per snapshot, respectively. We set $\delta = 0.05$ and $\epsilon = 0.5$ in Algorithm \ref{alg:ridge_leverage_score_CSS} to perform CSS. The reconstruction result of different time steps with target rank $k=50$ is shown in Figs.\ \ref{fig:channelflow_1}--\ref{fig:channelflow_3}. We find that most turbulence features are well captured in the reconstruction result. However, we also observe ripple shapes that appear in some reconstruction images (Step 220 in Fig.\ \ref{fig:channelflow_1}, Step 920 in Fig.~\ref{fig:channelflow_2}) causing a larger reconstruction error.

To evaluate the performance of our approach, we measure the effect of the target rank $k$ on the relative error (measured in terms of the Frobenius norm) of different coefficient computation algorithms as shown in Table \ref{tab:channelflow}. We select the truncated SVD method and the single-pass randomized SVD \cite{yu2017single} using different target ranks as the benchmarks for comparison. For the residual-based CSS (Algorithm \ref{alg:CSS_residual}), to ensure its memory usage of column basis is the same as other approaches, we only expand its column subset until $k$ columns are selected and any new column will not be added to the subset. The same treatment is also used in the remaining experiments of this paper. For each case, we also record the relative error estimated using the single-pass Hutch++ method. As expected, we find that the reconstruction accuracy keeps improving along with the target rank. Given the same rank, we use 
the color red to denote
the exact and estimated relative Frobenius error with the least value. We find that the accuracy difference between the coefficients obtained from the four strategies of Section \ref{subsec:coeff_compute} is very small in this experiment. In addition, the single-pass Hutch++ method accurately estimates the relative Frobenius error and determines the coefficient matrix with the least reconstruction error.

To isolate just the impact of the ridge leverage score-based CSS algorithm, we first compare it with the other three CSS methods based on CPQR, residual norm (Algorithm \ref{alg:CSS_residual}), and standard leverage score sampling \cite{jolliffe1972discarding} in a multi-pass setting. In these other approaches, we first obtain the indices of the selected columns $\mathcal{J}$ and then compute the coefficient matrix by directly solving Eq.~\eqref{eq:coeff_ID1} without sketching. The ridge leverage score-based CSS performs slightly worse than the standard leverage score method, which is reasonable since the score is approximated in a streaming setting. It also performs better than the residual-based CSS but the CPQR outperforms all the other approaches. In Table~\ref{tab:channelflow}, we then compare the ridge leverage score-based CSS (Algorithm \ref{alg:ridge_leverage_score_CSS}) with the residual-based CSS (Algorithm \ref{alg:CSS_residual}) in an online setting where the coefficient matrix is computed simultaneously with CSS using Algs. \ref{alg:coeff_update1}--\ref{alg:coeff_update4}. In most cases, our choice of CSS approach obtains better reconstruction accuracy than the residual-based CSS method.

There are a few interesting observations to make from Table~\ref{tab:channelflow}. First, comparing SVD to ID methods, there is little difference at low ranks but a larger difference at high ranks where the SVD method starts to outperform the ID vastly. For example, looking at the $256\times 256$ grid data, at rank 5, the ratio of CPQR to SVD error is $2.89/1.95\approx 1.48\times$, whereas at rank $100$ the ratio of error is $.67/.13\approx 5.15\times$. This suggests that the benefits of an ID method (e.g., interpretability) may be worth it for low ranks, but less so for higher ranks.  Second, looking now only among the ID methods, we also see a stronger gap emerge only at larger ranks. Again looking at the $256\times 256$ grid data, the ratio of rID to CPQR error is just $3.28/2.89\approx 1.13\times$ at rank 5 but worsens to $1.33/0.67\approx 1.99\times$ at rank 100. Overall, for ranks up to 10, the one-pass rID is within $2\times$ to $3\times$ the accuracy of the (multi-pass) SVD, suggesting that there is little room for improvement in this regime. As the rank increases, the rID error improves but only mildly, in comparison to the multi-pass SVD error which greatly improves.

\newcommand{\red}[1]{\textcolor{red}{#1}} 
\newcommand{\blue}[1]{\textcolor{blue}{#1}} 
\newcommand{\best}[1]{\textcolor{red}{#1}} 
\sisetup{
text-series-to-math = true,
propagate-math-font = true
}
\newcommand{\frmt}[1]{\num[exponent-mode = fixed, fixed-exponent = 0]{#1}\si{\percent}} 
\begin{sidewaystable}[!htbp]
        \small
        \caption{The relative error of randomized ID with different coefficient computation algorithms and different numbers of new columns. The estimated error using single-pass Hutch++ is included in the bracket. The testing is performed on the turbulence channel flow dataset}
        \label{tab:channelflow}
        \begin{adjustbox}{width = 0.9\textwidth, center}
                \centering
                \begin{tabular}{ccllllllll}
    \toprule
            &           &                                                   &                              & \multicolumn{6}{c}{Target rank}                                                                                                                                                                                                                                                      \\
    \cmidrule(l){5-10}
    Dataset & \# passes & Column selection                                  & Coeff.\ update               & 5                                          & 10                                         & 20                                         & 40                                          & 50                                          & 100                                        \\
    \midrule
    \multirow{14}{*}{\shortstack{Channel flow                                                                                                                                                                                                                                                                                                                                              \\$64 \times 64 \times 1000$}}
            & $2$       & \multicolumn{2}{c}{Truncated SVD}                 & \frmt{1.61}                  & \frmt{1.17}                                & \frmt{0.78}                                & \frmt{0.41}                                & \frmt{0.31}                                 & \frmt{0.11}                                                                              \\
            & $1$       & \multicolumn{2}{c}{Randomized SVD}                & \frmt{1.71}                  & \frmt{1.40}                                & \frmt{1.21}                                & \frmt{0.72}                                & \frmt{0.58}                                 & \frmt{0.23}                                                                              \\
    \cmidrule(l){2-10}
            & $k$       & CPQR                                        & Least squares (Eq.~\eqref{eq:coeff_ID2})  & \frmt{2.51}                                & \frmt{2.30}                                & \frmt{1.58}                                & \frmt{1.20}                                 & \frmt{0.97}                                 & \frmt{0.47}                                \\
            & $2$       & Leverage score                                    & Least squares (Eq.~\eqref{eq:coeff_ID2})  & \frmt{3.17}                                & \frmt{2.84}                                & \frmt{1.86}                                & \frmt{1.94}                                 & \frmt{1.22}                                 & \frmt{1.15}                                \\
            & $2$       & Residual-based (Algorithm \ref{alg:CSS_residual})                                & Least squares (Eq.~\eqref{eq:coeff_ID2})  & \frmt{3.47}                                & \frmt{3.05}                                & \frmt{2.75}                                & \frmt{2.61}                                 & \frmt{1.39}                                 & \frmt{1.22}                                \\
            & $2$       & Ridge leverage score-based (Algorithm \ref{alg:rID_overview})                     & Least squares (Eq.~\eqref{eq:coeff_ID2})  & \frmt{3.25}                                & \frmt{3.01}                                & \frmt{2.12}                                & \frmt{2.04}                                 & \frmt{1.35}                                 & \frmt{1.21}                                \\
    \cmidrule(l){2-10}
            & $1$       & Residual-based (Algorithm \ref{alg:CSS_residual})                                & Algorithm \ref{alg:coeff_update1}              & \frmt{4.45}(\frmt{4.67})                   & \frmt{3.35}(\frmt{3.97})                   & \best{\frmt{3.78}}(\best{\frmt{4.25}}) & \best{\frmt{3.62}}(\best{\frmt{4.35}}) & \best{\frmt{2.75}}(\best{\frmt{3.94}}) & \frmt{2.29}(\frmt{2.82})                   \\
            & $1$       & Residual-based (Algorithm \ref{alg:CSS_residual})                                & Algorithm \ref{alg:coeff_update2} & \best{\frmt{4.15}}(\best{\frmt{4.34}}) & \best{\frmt{3.17}}(\best{\frmt{3.52}}) & \frmt{4.75}(\frmt{4.96})                   & \frmt{5.73}(\frmt{6.34})                    & \frmt{6.91}(\frmt{7.83})                    & \frmt{4.73}(\frmt{5.02})                   \\
            & $1$       & Residual-based (Algorithm \ref{alg:CSS_residual})                                & Algorithm \ref{alg:coeff_update3}     & \frmt{4.69}(\frmt{4.82})                   & \frmt{4.52}(\frmt{4.71})                   & \frmt{4.18}(\frmt{4.53})                   & \frmt{5.78}(\frmt{6.26})                    & \frmt{5.81}(\frmt{7.03})                    & \frmt{3.35}(\frmt{3.67})                   \\
            & $1$       & Residual-based (Algorithm \ref{alg:CSS_residual})                                & Algorithm \ref{alg:coeff_update4}                    & \frmt{4.33}(\frmt{4.57})                   & \frmt{3.95}(\frmt{4.31})                   & \frmt{4.01}(\frmt{4.71})                   & \frmt{4.53}(\frmt{5.18})                    & \frmt{3.78}(\frmt{4.16})                    & \best{\frmt{2.15}}(\best{\frmt{2.78}}) \\
    \cmidrule(l){2-10}
            & $1$       & Ridge leverage score-based (Algorithm \ref{alg:rID_overview})                     & Algorithm \ref{alg:coeff_update1}              & \frmt{3.91}(\frmt{4.06})                   & \frmt{3.51}(\frmt{3.47})                   & \frmt{3.27}(\frmt{3.25})                   & \frmt{3.24} (\frmt{3.31})                   & \frmt{2.35}(\frmt{2.43})                    & \frmt{2.11}(\frmt{2.17})                   \\
            & $1$       & Ridge leverage score-based (Algorithm \ref{alg:rID_overview})                     & Algorithm \ref{alg:coeff_update2} & \best{\frmt{3.83}}(\best{\frmt{3.99}}) & \frmt{3.43}(\frmt{3.52})                   & \frmt{3.57}(\frmt{3.93})                   & \frmt{3.92}(\frmt{4.02})                    & \frmt{3.19} (\frmt{3.24})                   & \frmt{2.55}(\frmt{2.42})                   \\
            & $1$       & Ridge leverage score-based (Algorithm \ref{alg:rID_overview})                     & Algorithm \ref{alg:coeff_update3}     & \frmt{4.17}(\frmt{4.26})                   & \best{\frmt{3.15}}(\best{\frmt{3.22}}) & \best{\frmt{3.21}}(\best{\frmt{3.12}}) & \best{\frmt{3.13}}(\best{\frmt{3.22}})  & \best{\frmt{1.95}}(\best{\frmt{2.01}})  & \best{\frmt{1.96}}(\best{\frmt{2.07}}) \\
            & $1$       & Ridge leverage score-based (Algorithm \ref{alg:rID_overview})                     & Algorithm \ref{alg:coeff_update4}                    & \frmt{4.25}(\frmt{4.35})                   & \frmt{3.49}(\frmt{3.58})                   & \frmt{3.81}(\frmt{3.74})                   & \frmt{3.75}(\frmt{3.87})                    & \frmt{3.20} (\frmt{3.37})                   & \frmt{2.84}(\frmt{3.11})                   \\
    \cmidrule(l){2-10}
            &           & \multicolumn{2}{c}{Update decision}               & Algorithm \ref{alg:coeff_update2} & Algorithm \ref{alg:coeff_update3}               & Algorithm \ref{alg:coeff_update3}                            & Algorithm \ref{alg:coeff_update3}                            & Algorithm \ref{alg:coeff_update3}                             & Algorithm \ref{alg:coeff_update3}                                                                                \\
            &           & \multicolumn{2}{c}{Select best update every time} & \frmt{4.02}                  & \frmt{3.54}                                & \frmt{3.42}                                & \frmt{3.16}                                & \frmt{2.84}                                 & \frmt{1.96}                                                                              \\
    \midrule
    \multirow{14}{*}{\shortstack{Channel flow                                                                                                                                                                                                                                                                                                                                              \\$128 \times 128 \times 1000$}}
            & $2$       & \multicolumn{2}{c}{Truncated SVD}                 & \frmt{1.76}                  & \frmt{1.34}                                & \frmt{0.91}                                & \frmt{0.47}                                & \frmt{0.34}                                 & \frmt{0.10}                                                                              \\
            & $1$       & \multicolumn{2}{c}{Randomized SVD}                & \frmt{2.00}                  & \frmt{1.58}                                & \frmt{1.41}                                & \frmt{0.83}                                & \frmt{0.65}                                 & \frmt{0.23}                                                                              \\
    \cmidrule(l){2-10}
            & $k$       & CPQR                                        & Least squares (Eq.~\eqref{eq:coeff_ID2})  & \frmt{2.67}                                & \frmt{2.25}                                & \frmt{2.01}                                & \frmt{1.28}                                 & \frmt{1.04}                                 & \frmt{0.49}                                \\
            & $2$       & Leverage score                                    & Least squares (Eq.~\eqref{eq:coeff_ID2})  & \frmt{2.89}                                & \frmt{2.30}                                & \frmt{2.08}                                & \frmt{1.36}                                 & \frmt{1.17}                                 & \frmt{0.47}                                \\
            & $2$       & Residual-based (Algorithm \ref{alg:CSS_residual})                                & Least squares (Eq.~\eqref{eq:coeff_ID2})  & \frmt{3.12}                                & \frmt{2.55}                                & \frmt{2.13}                                & \frmt{1.46}                                 & \frmt{1.32}                                 & \frmt{1.18}                                \\
            & $2$       & Ridge leverage score-based (Algorithm \ref{alg:rID_overview})                     & Least squares (Eq.~\eqref{eq:coeff_ID2})  & \frmt{3.05}                                & \frmt{2.50}                                & \frmt{2.21}                                & \frmt{1.44}                                 & \frmt{1.22}                                 & \frmt{1.15}                                \\
    \cmidrule(l){2-10}
            & $1$       & Residual-based (Algorithm \ref{alg:CSS_residual})                                & Algorithm \ref{alg:coeff_update1}              & \frmt{3.62}(\frmt{4.61})                   & \frmt{3.98}(\frmt{4.82})                   & \frmt{3.52}(\frmt{4.58})                   & \frmt{4.12}(\frmt{4.36})                    & \frmt{4.54}(\frmt{5.08})                    & \frmt{2.84}(\frmt{3.95})                   \\
            & $1$       & Residual-based (Algorithm \ref{alg:CSS_residual})                                & Algorithm \ref{alg:coeff_update2} & \best{\frmt{3.50}}(\best{\frmt{3.92}}) & \frmt{5.05}(\frmt{5.72})                   & \best{\frmt{3.27}}(\best{\frmt{3.82}}) & \frmt{3.82}(\frmt{4.81})                    & \frmt{3.25}(\frmt{3.35})                    & \frmt{3.54}(\frmt{4.24})                   \\
            & $1$       & Residual-based (Algorithm \ref{alg:CSS_residual})                                & Algorithm \ref{alg:coeff_update3}     & \frmt{4.15}(\frmt{4.53})                   & \frmt{4.18}(\frmt{4.61})                   & \frmt{3.71}(\frmt{4.35})                   & \best{\frmt{3.73}}(\best{\frmt{4.72}})  & \frmt{2.65}(\frmt{3.80})                    & \best{\frmt{2.88}}(\best{\frmt{3.54}}) \\
            & $1$       & Residual-based (Algorithm \ref{alg:CSS_residual})                                & Algorithm \ref{alg:coeff_update4}                    & \frmt{3.98}(\frmt{4.29})                   & \best{\frmt{4.24}}(\best{\frmt{4.52}}) & \frmt{4.05}(\frmt{4.76})                   & \frmt{5.28}(\frmt{6.10})                    & \best{\frmt{2.41}}(\best{\frmt{3.22}})  & \frmt{3.18}(\frmt{3.85})                   \\
    \cmidrule(l){2-10}
            & $1$       & Ridge leverage score-based (Algorithm \ref{alg:rID_overview})                     & Algorithm \ref{alg:coeff_update1}              & \frmt{3.84}(\frmt{3.96})                   & \frmt{3.33}(\frmt{3.56})                   & \frmt{2.97}(\frmt{3.02})                   & \frmt{3.31}(\frmt{3.42})                    & \frmt{3.48}(\frmt{3.57})                    & \frmt{2.33}(\frmt{2.41})                   \\
            & $1$       & Ridge leverage score-based (Algorithm \ref{alg:rID_overview})                     & Algorithm \ref{alg:coeff_update2} & \best{\frmt{3.21}}(\best{\frmt{3.33}}) & \frmt{3.94}(\frmt{3.76})                   & \best{\frmt{2.79}}(\best{\frmt{3.02}}) & \frmt{3.35}(\frmt{3.67})                    & \frmt{2.85}(\frmt{2.67})                    & \frmt{2.46}(\frmt{3.04})                   \\
            & $1$       & Ridge leverage score-based (Algorithm \ref{alg:rID_overview})                     & Algorithm \ref{alg:coeff_update3}     & \frmt{3.72}(\frmt{3.98})                   & \frmt{3.18}(\frmt{3.38})                   & \frmt{3.03}(\frmt{3.65})                   & \best{\frmt{3.02}}(\best{\frmt{4.19}})  & \frmt{3.00}(\frmt{3.12})                    & \best{\frmt{1.76}}(\best{\frmt{2.06}}) \\
            & $1$       & Ridge leverage score-based (Algorithm \ref{alg:rID_overview})                     & Algorithm \ref{alg:coeff_update4}                    & \frmt{3.46}(\frmt{3.76})                   & \best{\frmt{2.96}}(\best{\frmt{3.10}}) & \frmt{2.85}(\frmt{3.16})                   & \frmt{3.11}(\frmt{3.72})                    & \best{\frmt{1.99}}(\best{\frmt{2.44}})  & \frmt{1.98}(\frmt{2.35})                   \\
    \cmidrule(l){2-10}
            &           & \multicolumn{2}{c}{Update decision}               & Algorithm \ref{alg:coeff_update2} & Algorithm \ref{alg:coeff_update4}                                & Algorithm \ref{alg:coeff_update2}               & Algorithm \ref{alg:coeff_update3}                                        & Algorithm \ref{alg:coeff_update4}                                   & Algorithm \ref{alg:coeff_update3}                                                                                      \\
            &           & \multicolumn{2}{c}{Select best update every time} & \frmt{3.85}                  & \frmt{2.68}                                & \frmt{3.15}                                & \frmt{2.77}                                & \frmt{2.42}                                 & \frmt{1.93}                                                                              \\
    \midrule
    \multirow{14}{*}{\shortstack{Channel flow                                                                                                                                                                                                                                                                                                                                              \\$256 \times 256 \times 1000$}}
            & $2$       & \multicolumn{2}{c}{Truncated SVD}                 & \frmt{1.95}                  & \frmt{1.37}                                & \frmt{1.01}                                & \frmt{0.52}                                & \frmt{0.39}                                 & \frmt{0.13}                                                                              \\
            & $1$       & \multicolumn{2}{c}{Randomized SVD}                & \frmt{2.18}                  & \frmt{1.82}                                & \frmt{1.62}                                & \frmt{0.96}                                & \frmt{0.78}                                 & \frmt{0.29}                                                                              \\
    \cmidrule(l){2-10}
            & $k$       & CPQR                                        & Least squares (Eq.~\eqref{eq:coeff_ID2})  & \frmt{2.89}                                & \frmt{1.76}                                & \frmt{1.55}                                & \frmt{1.19}                                 & \frmt{0.98}                                 & \frmt{0.67}                                \\
            & $2$       & Leverage score                                    & Least squares (Eq.~\eqref{eq:coeff_ID2})  & \frmt{2.84}                                & \frmt{1.89}                                & \frmt{1.69}                                & \frmt{1.37}                                 & \frmt{0.97}                                 & \frmt{0.79}                                \\
            & $2$       & Residual-based (Algorithm \ref{alg:CSS_residual})                                & Least squares (Eq.~\eqref{eq:coeff_ID2})  & \frmt{2.95}                                & \frmt{2.17}                                & \frmt{1.74}                                & \frmt{1.43}                                 & \frmt{1.15}                                 & \frmt{1.27}                                \\
            & $2$       & Ridge leverage score-based (Algorithm \ref{alg:rID_overview})                     & Least squares (Eq.~\eqref{eq:coeff_ID2})  & \frmt{2.90}                                & \frmt{1.95}                                & \frmt{1.80}                                & \frmt{1.39}                                 & \frmt{1.02}                                 & \frmt{0.94}                                \\
    \cmidrule(l){2-10}
            & $1$       & Residual-based (Algorithm \ref{alg:CSS_residual})                                & Algorithm \ref{alg:coeff_update1}              & \frmt{4.14}(\frmt{4.35})                   & \frmt{3.54}(\frmt{3.72})                   & \frmt{2.90}(\frmt{3.71})                   & \frmt{2.25}(\frmt{2.76})                    & \frmt{4.21}(\frmt{4.67})                    & \frmt{2.97}(\frmt{4.75})                   \\
            & $1$       & Residual-based (Algorithm \ref{alg:CSS_residual})                                & Algorithm \ref{alg:coeff_update2} & \frmt{3.76}(\frmt{4.05})                   & \best{\frmt{2.86}}(\best{\frmt{3.50}}) & \frmt{3.19}(\frmt{3.92})                   & \frmt{2.70}(\frmt{3.77})                    & \frmt{3.14}(\frmt{4.18})                    & \frmt{3.51}(\frmt{4.61})                   \\
            & $1$       & Residual-based (Algorithm \ref{alg:CSS_residual})                                & Algorithm \ref{alg:coeff_update3}     & \frmt{3.98}(\frmt{4.27})                   & \frmt{3.12}(\frmt{3.79})                   & \best{\frmt{2.76}}(\best{\frmt{3.74}}) & \frmt{2.40}(\frmt{3.24})                    & \best{\frmt{2.06}}(\best{\frmt{2.69}})  & \frmt{3.11}(\frmt{3.52})                   \\
            & $1$       & Residual-based (Algorithm \ref{alg:CSS_residual})                                & Algorithm \ref{alg:coeff_update4}                    & \best{\frmt{3.74}}(\best{\frmt{4.18}}) & \frmt{3.57}(\frmt{4.15})                   & \frmt{2.85}(\frmt{4.25})                   & \best{\frmt{1.54}}(\best{\frmt{1.93}})  & \frmt{3.10}(\frmt{4.37})                    & \best{\frmt{2.51}}(\best{\frmt{3.17}}) \\
    \cmidrule(l){2-10}
            & $1$       & Ridge leverage score-based (Algorithm \ref{alg:rID_overview})                     & Algorithm \ref{alg:coeff_update1}              & \frmt{3.06}(\frmt{3.25})                   & \frmt{2.32}(\frmt{2.49})                   & \frmt{2.85}(\frmt{2.93})                   & \frmt{2.54}(\frmt{2.63})                    & \frmt{2.63}(\frmt{2.74})                    & \frmt{2.39}(\frmt{2.48})                   \\
            & $1$       & Ridge leverage score-based (Algorithm \ref{alg:rID_overview})                     & Algorithm \ref{alg:coeff_update2} & \frmt{3.26}(\frmt{3.54})                   & \best{\frmt{2.18}}(\best{\frmt{2.36}}) & \frmt{3.19}(\frmt{3.24})                   & \frmt{2.17}(\frmt{2.07})                    & \frmt{2.64}(\frmt{3.38})                    & \frmt{2.35}(\frmt{2.61})                   \\
            & $1$       & Ridge leverage score-based (Algorithm \ref{alg:rID_overview})                     & Algorithm \ref{alg:coeff_update3}     & \frmt{3.38}(\frmt{3.72})                   & \frmt{2.92}(\frmt{2.74})                   & \best{\frmt{2.76}}(\best{\frmt{2.71}}) & \frmt{2.24}(\frmt{2.64})                    & \best{\frmt{2.76}}(\best{\frmt{2.84}})  & \frmt{2.71}(\frmt{2.59})                   \\
            & $1$       & Ridge leverage score-based (Algorithm \ref{alg:rID_overview})                     & Algorithm \ref{alg:coeff_update4}                    & \best{\frmt{3.04}}(\best{\frmt{3.17}}) & \frmt{2.67}(\frmt{2.13})                   & \frmt{2.85}(\frmt{3.25})                   & \best{\frmt{2.13}}(\best{\frmt{1.96}})  & \frmt{2.92}(\frmt{3.27})                    & \best{\frmt{1.51}}(\best{\frmt{1.98}}) \\
    \cmidrule(l){2-10}
            &           & \multicolumn{2}{c}{Update decision}               & Algorithm \ref{alg:coeff_update4}                    & Algorithm \ref{alg:coeff_update2}               & Algorithm \ref{alg:coeff_update3}                                   & Algorithm \ref{alg:coeff_update4}                                  & Algorithm \ref{alg:coeff_update3}                                        & Algorithm \ref{alg:coeff_update4}                                                                                     \\
            &           & \multicolumn{2}{c}{Select best update every time} & \frmt{3.28}                  & \frmt{3.01}                                & \frmt{2.16}                                & \frmt{1.75}                                & \frmt{1.42}                                 & \frmt{1.33}                                                                              \\
    \bottomrule
\end{tabular}
        \end{adjustbox}
\end{sidewaystable}

\subsubsection{Incorporating gradient information}
\begin{figure}[!htbp]
        \centering
        \includegraphics[width=\textwidth]{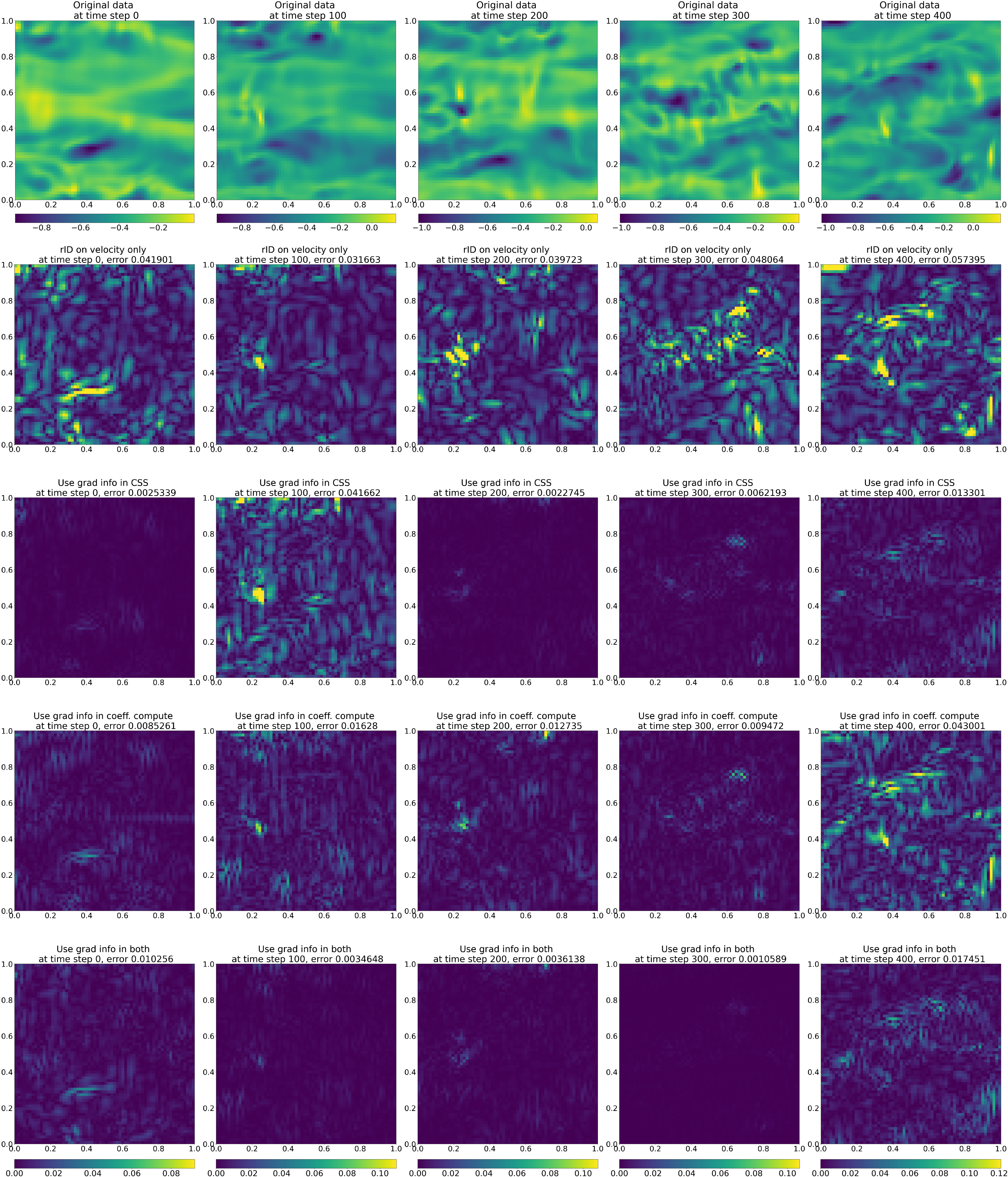}
        \caption{Rank $k = 50$, the $z$-direction vorticity reconstruction of the turbulence flow data over a $64 \times 64$ grid at different time steps. The first row shows the original data. The remaining 4 rows, from top to bottom, show the nodal absolute reconstruction error using: randomized ID on velocity only; randomized ID with gradient information in CSS; randomized ID with gradient information in coefficient computation; and randomized ID with velocity and estimated gradient information.}
        \label{fig:channelflow_1_vortocity}
   \end{figure}

For the turbulence channel flow dataset, it is also important to maintain an accurate reconstruction of vorticity. Since the computation of vorticity utilizes the gradient of the velocity, we need to account for the gradient reconstruction and thus compared four testing cases: (1) no gradient estimation; (2) using gradient estimation in coefficient computation; (3) using gradient estimation in CSS; and (4) using gradient estimation in both CSS and coefficient computation. 
In Fig.~\ref{fig:channelflow_1_vortocity}, we plot the original $z$-direction vorticity ($1^\text{st}$ row) and compare the reconstructed vorticity field using the aforementioned four cases at different time steps. To highlight the difference among these four cases, we compute and plot the absolute reconstruction error at each node of the grid for comparison ($2^\text{nd}$--$5^\text{th}$ rows). Compared to the $2^\text{nd}$ row, the vorticity reconstruction accuracy is obviously improved when the gradient is used during either CSS or coefficient computation ($3^\text{rd}$--$5^\text{th}$ rows). When the gradient is used only in CSS ($3^\text{rd}$ row), the column basis is improved with more representative columns and thus the reconstruction error at time step $0$ and $200$ are minor since they are selected as column basis. When the gradient is only used for coefficient computation ($4^\text{th}$ row), we find that the reconstruction accuracy is increased due to better coefficients. These factors contribute to the significant improvement of the vorticity reconstruction when the gradient information is used in both CSS and coefficient computation procedures ($5^\text{th}$ row). 

In Table \ref{tab:channelflow2}, we summarize the reconstruction error of both the velocity and vorticity fields. We find that the vorticity accuracy improves when the gradient information is used in any fashion (cases 2--4), as to be expected, and is typically best when the gradient is used as much as possible (case 4), followed by case 2, case 3 and worst in case 1. 
Interestingly, including the gradient information does not necessarily hurt the velocity reconstruction as one might expect. That is, case 3 regularly improves on case 1 in terms of velocity accuracy. We do not have a full explanation, but it it could be related to the fact that the column selection of rID is not optimal and that by adding in the gradient term adds information. It could also be possible that case 4 would improve even further but is being limited by the sketching. Overall, the results suggest using case 3 if one only cares about velocity accuracy, and using case 4 if one cares about both velocity and vorticity.

\begin{table}[!htbp]
        \resizebox{\textwidth}{!}{
                \centering
                \begin{tabular}{cllllllll}
                        \toprule
                                &                             &                          & \multicolumn{6}{c}{Target rank}                                                                                                                                                      \\
                        \cmidrule(l){4-9}
                        Dataset & Column selection            & Coeff.\ update  criteria & 5                           & 10                          & 20                         & 40                         & 50                         & 100                        \\
                        \midrule
\multirow{4}{*}{\shortstack{Channel flow   \\$64 \times 64 \times 1000$}}
& Only velocity term          & Only velocity vector     & \frmt{4.02} \& \frmt{11.08} & \frmt{3.54} \& \frmt{10.93} & \frmt{3.42} \& \frmt{9.64} & \frmt{3.16} \& \frmt{9.33} & \frmt{2.84} \& \frmt{9.27} & \frmt{1.96} \& \frmt{7.42} \\
& Only velocity term          & With velocity gradient   & \frmt{4.17} \& \frmt{7.83}  & \frmt{3.95} \& \frmt{6.93}  & \frmt{3.87} \& \frmt{6.12} & \frmt{3.24} \& \best{\frmt{5.27}} & \frmt{3.06} \& \frmt{5.19} & \frmt{2.21} \& \frmt{4.02} \\
& With velocity gradient term & Only velocity vector     & \frmt{3.85} \& \frmt{9.32}  & \best{\frmt{3.32}} \& \frmt{9.15}  & \best{\frmt{2.97}} \& \frmt{8.37} & \best{\frmt{2.85}} \& \frmt{7.94} & \best{\frmt{2.72}} \& \frmt{7.32} & \best{\frmt{1.83}} \& \frmt{6.41} \\
& With velocity gradient term & With velocity gradient   & \best{\frmt{3.25}} \& \best{\frmt{6.25}}  & \frmt{3.77} \& \best{\frmt{6.07}}  & \frmt{3.64} \& \best{\frmt{5.91}} & \frmt{3.19} \& \frmt{5.42} & \frmt{2.98} \& \best{\frmt{5.01}} & \frmt{2.07} \& \best{\frmt{3.74}} \\
                        \midrule
                        \multirow{4}{*}{\shortstack{Channel flow \\$128 \times 128 \times 1000$}}
                                & Only velocity term          & Only velocity vector     & \frmt{3.85} \& \frmt{10.26} & \frmt{3.68} \& \frmt{9.92}  & \frmt{3.15} \& \frmt{8.17} & \best{\frmt{2.77}} \& \frmt{7.02} & \frmt{2.42} \& \frmt{6.94} & \frmt{1.93} \& \frmt{6.81} \\
                                & Only velocity term          & With velocity gradient   & \frmt{4.39} \& \best{\frmt{5.17}}  & \frmt{4.15} \& \best{\frmt{5.01}}  & \frmt{3.87} \& \best{\frmt{4.93}} & \frmt{3.34} \& \frmt{4.77} & \frmt{2.94} \& \frmt{4.67} & \frmt{2.19} \& \frmt{4.13} \\
                                & With velocity gradient term & Only velocity vector     & \frmt{3.42} \& \frmt{7.85}  & \best{\frmt{3.17}} \& \frmt{6.94}  & \best{\frmt{2.94}} \& \frmt{6.35} & \frmt{2.67} \& \frmt{6.01} & \best{\frmt{2.13}} \& \frmt{5.85} & \best{\frmt{1.78}} \& \frmt{5.14} \\
                                & With velocity gradient term & With velocity gradient   & \best{\frmt{3.23}} \& \frmt{5.31}  & \frmt{3.35} \& \frmt{5.11}  & \frmt{3.06} \& \frmt{5.02} & \frmt{2.93} \& \best{\frmt{4.95}} & \frmt{2.57} \& \best{\frmt{4.47}} & \frmt{2.08} \& \best{\frmt{3.92}} \\
                        \midrule
                        \multirow{4}{*}{\shortstack{Channel flow  \\$256 \times 256 \times 1000$}}
                                & Only velocity term          & Only velocity vector     & \frmt{3.28} \& \frmt{9.02}  & \frmt{3.01} \& \frmt{8.66}  & \frmt{2.16} \& \frmt{8.12} & \frmt{1.75} \& \frmt{7.02} & \frmt{1.42} \& \frmt{6.98} & \frmt{1.33} \& \frmt{6.33} \\
                                & Only velocity term          & With velocity gradient   & \frmt{3.84} \& \frmt{5.72}  & \frmt{3.29} \& \frmt{5.17}  & \frmt{2.34} \& \frmt{5.01} & \frmt{2.08} \& \frmt{4.95} & \frmt{1.87} \& \frmt{4.72} & \frmt{1.64} \& \frmt{4.17} \\
                                & With velocity gradient term & Only velocity vector     & \frmt{3.15} \& \frmt{7.15}  & \best{\frmt{2.87}} \& \frmt{6.83}  & \best{\frmt{2.03}} \& \frmt{6.35} & \best{\frmt{1.65}} \& \frmt{5.80} & \best{\frmt{1.22}} \& \frmt{5.67} & \best{\frmt{1.17}} \& \frmt{5.35} \\
                                & With velocity gradient term & With velocity gradient   & \best{\frmt{2.90}} \& \best{\frmt{5.06}}  & \frmt{2.93} \& \best{\frmt{4.89}}  & \frmt{2.12} \& \best{\frmt{4.75}} & \frmt{1.89} \& \best{\frmt{4.23}} & \frmt{1.45} \& \frmt{4.13} & \frmt{1.48} \& \best{\frmt{4.01}} \\
                        \bottomrule
                \end{tabular}
        }
        \caption{The relative error of velocity ($1^\text{st}$ value) and vorticity ($2^\text{nd}$ value) with different optimization functions and CSS criteria. The testing is performed on the turbulence channel flow dataset.}
        \label{tab:channelflow2}
\end{table}

\subsection{Ignition dataset}
\label{subsec:ignition}
\begin{figure}[!ht]
        \centering
        \includegraphics[width =\textwidth]{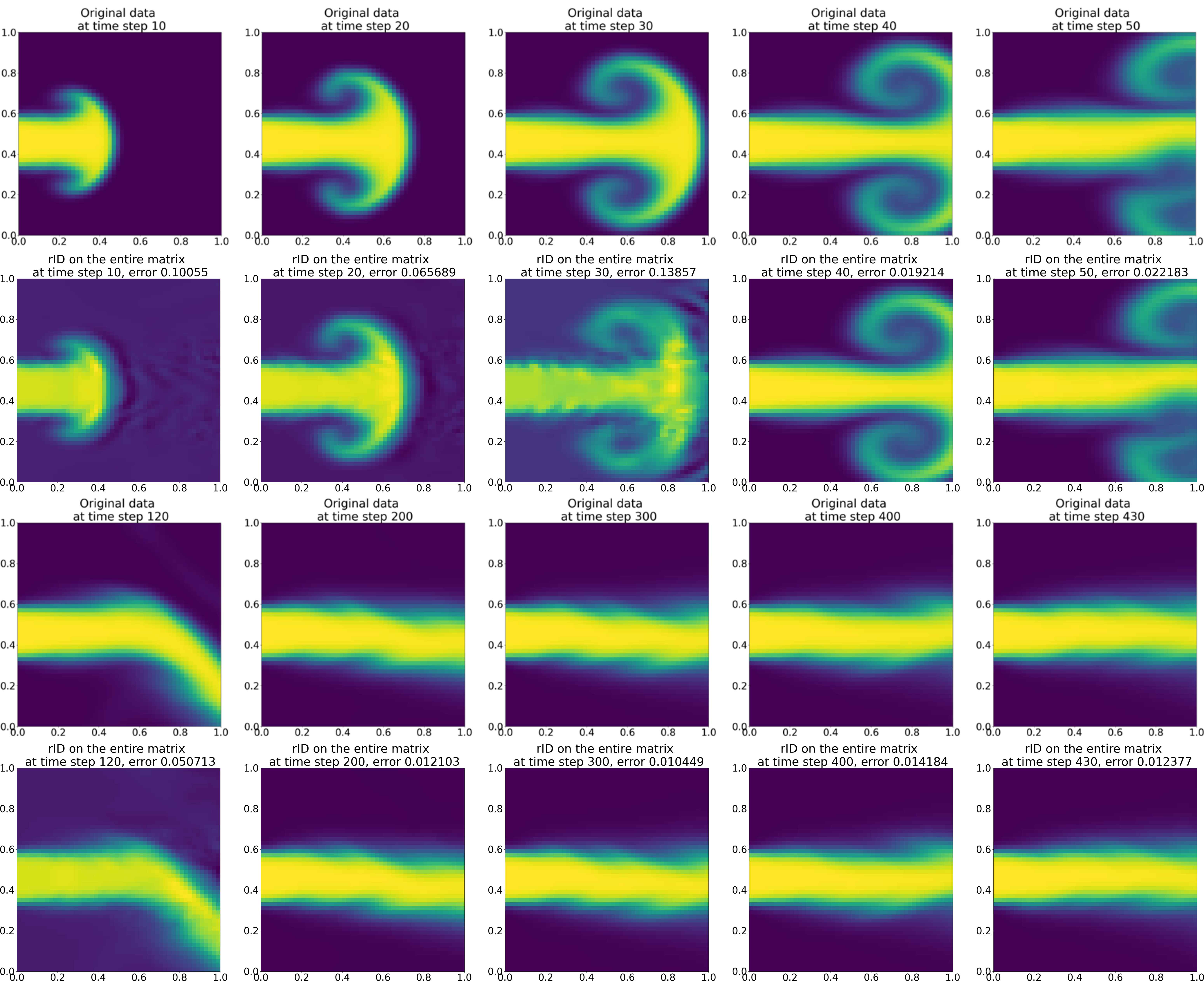}
        \caption{Rank $k=40$ reconstruction of ignition data on a grid of size $50 \times 50$ at different time steps. For each two rows, the top one represents the original data while the bottom one represents the reconstruction data.}
        \label{fig:ignition_1}
\end{figure}

In our second experiment, we evaluate the performance of our algorithm using the ignition data obtained from a jet diffusion simulation on a uniform grid of size $50 \times 50$. In particular, the dataset is convection-dominated for the early ignition stage until the jet flame reaches a steady state. Since the flame profile changes drastically during ignition, the steady-state data cannot be well reconstructed with only the early ignition stage data, or vice versa. Therefore, the approximation accuracy is sensitive to the selection of column basis during the streaming compression and it is essential for the CSS algorithm to select the column basis containing the data in both ignition and steady-state stages. The reconstruction at several different time points is shown in Fig. \ref{fig:ignition_1}.
\begin{table}[!htb]
        \resizebox{\textwidth}{!}{
                \centering
                \begin{tabular}{ccllllllll}
        \toprule
                &           &                                                   &                              & \multicolumn{6}{c}{Target rank}                                                                                                                                                                                                                                                              \\
        \cmidrule(l){5-10}
        Dataset & \# passes & Column selection                                  & Coeff.\ update               & 5                                            & 10                                           & 20                                           & 40                                           & 50                                           & 100                                        \\
        \midrule
        \multirow{14}{*}{\shortstack{Ignition grid \\ $50 \times 50 \times 450$}}
                & $2$       & \multicolumn{2}{c}{Truncated SVD}                 & \frmt{9.63}                  & \frmt{5.99}                                  & \frmt{3.06}                                  & \frmt{1.21}                                  & \frmt{0.83}                                  & \frmt{0.20}                                                                               \\
                & $1$       & \multicolumn{2}{c}{Randomized SVD}                & \frmt{17.7}                  & \frmt{11.9}                                  & \frmt{8.90}                                  & \frmt{4.00}                                  & \frmt{3.09}                                  & \frmt{0.82}                                                                               \\
        \cmidrule(l){2-10}
                & $k$       & CPQR                                        & Least squares (Eq.~\eqref{eq:coeff_ID2})  & \frmt{17.45}                                 & \frmt{11.75}                                 & \frmt{7.85}                                  & \frmt{3.54}                                  & \frmt{2.52}                                  & \frmt{1.11}                                \\
                & $2$       & Leverage score                                    & Least squares (Eq.~\eqref{eq:coeff_ID2})  & \frmt{19.95}                                 & \frmt{14.64}                                 & \frmt{8.94}                                  & \frmt{4.19}                                  & \frmt{3.13}                                  & \frmt{2.91}                                \\
                & $2$       & Residual-based (Algorithm \ref{alg:CSS_residual})                                & Least squares (Eq.~\eqref{eq:coeff_ID2})  & \frmt{35.12}                                 & \frmt{28.73}                                 & \frmt{29.65}                                 & \frmt{18.34}                                 & \frmt{15.43}                                 & \frmt{6.19}                                \\
                & $2$       & Ridge leverage score-based (Algorithm \ref{alg:rID_overview})                     & Least squares (Eq.~\eqref{eq:coeff_ID2})  & \frmt{20.34}                                 & \frmt{16.39}                                 & \frmt{10.74}                                 & \frmt{7.38}                                  & \frmt{4.16}                                  & \frmt{3.70}                                \\
        \cmidrule(l){2-10}
                & $1$       & Residual-based (Algorithm \ref{alg:CSS_residual})                                & Algorithm \ref{alg:coeff_update1}              & \frmt{37.48}(\frmt{41.23})                   & \frmt{30.92}(\frmt{32.93})                   & \frmt{32.94}(\frmt{34.24})                   & \frmt{21.20}(\frmt{23.39})                   & \frmt{17.26}(\frmt{20.64})                   & \frmt{7.98}(\frmt{10.72})                  \\
                & $1$       & Residual-based (Algorithm \ref{alg:CSS_residual})                                & Algorithm \ref{alg:coeff_update2} & \frmt{38.10}(\frmt{40.25})                   & \best{\frmt{32.70}}(\best{\frmt{32.16}}) & \frmt{33.26}(\frmt{35.46})                   & \frmt{22.51}(\frmt{24.26})                   & \best{\frmt{16.98}}(\best{\frmt{19.34}}) & \frmt{8.27}(\frmt{10.54})                  \\
                & $1$       & Residual-based (Algorithm \ref{alg:CSS_residual})                                & Algorithm \ref{alg:coeff_update3}     & \frmt{37.24}(\frmt{39.76})                   & \frmt{31.83}(\frmt{33.72})                   & \frmt{33.19}(\frmt{36.75})                   & \best{\frmt{20.94}}(\best{\frmt{22.84}}) & \frmt{18.72}(\frmt{23.87})                   & \best{\frmt{7.15}}(\best{\frmt{9.26}}) \\
                & $1$       & Residual-based (Algorithm \ref{alg:CSS_residual})                                & Algorithm \ref{alg:coeff_update4}                    & \best{\frmt{36.30}}(\best{\frmt{39.69}}) & \frmt{30.80}(\frmt{33.25})                   & \best{\frmt{32.09}}(\best{\frmt{33.34}}) & \frmt{22.64}(\frmt{23.54})                   & \frmt{19.05}(\frmt{22.34})                   & \frmt{8.93}(\frmt{10.39})                  \\
        \cmidrule(l){2-10}
                & $1$       & Ridge leverage score-based (Algorithm \ref{alg:rID_overview})                     & Algorithm \ref{alg:coeff_update1}              & \frmt{25.26}(\frmt{27.13})                   & \frmt{21.81}(\frmt{22.15})                   & \frmt{13.36}(\frmt{14.97})                   & \best{\frmt{9.32}}(\best{\frmt{10.16}})                    & \frmt{5.92}(\frmt{6.03})                     & \frmt{4.85}(\frmt{5.03})                   \\
                & $1$       & Ridge leverage score-based (Algorithm \ref{alg:rID_overview})                     & Algorithm \ref{alg:coeff_update2} & \frmt{28.10}(\frmt{29.54})                   & \frmt{23.69}(\frmt{26.37})                   & \frmt{13.76}(\frmt{14.86})                   & \frmt{10.51}(\frmt{11.37})                   & \best{\frmt{4.98}}(\best{\frmt{5.34}})   & \frmt{3.91}(\frmt{4.28})                   \\
                & $1$       & Ridge leverage score-based (Algorithm \ref{alg:rID_overview})                     & Algorithm \ref{alg:coeff_update3}     & \frmt{27.16}(\frmt{28.93})                   & \frmt{24.17}(\frmt{26.92})                   & \frmt{13.79}(\frmt{13.85})                   & \frmt{11.98}(\frmt{12.85})                   & \frmt{5.17}(\frmt{7.18})                     & \frmt{4.93}(\frmt{5.44})                   \\
                & $1$       & Ridge leverage score-based (Algorithm \ref{alg:rID_overview})                     & Algorithm \ref{alg:coeff_update4}                    & \best{\frmt{23.19}}(\best{\frmt{24.73}}) & \best{\frmt{20.84}}(\best{\frmt{21.42}}) & \best{\frmt{12.40}}(\best{\frmt{13.54}}) & \frmt{10.64}(\frmt{12.49})                   & \frmt{5.35}(\frmt{6.24})                     & \best{\frmt{3.75}}(\best{\frmt{3.76}}) \\
        \cmidrule(l){2-10}
                &           & \multicolumn{2}{c}{Update decision}               & Algorithm \ref{alg:coeff_update4}                    & Algorithm \ref{alg:coeff_update4}                 & Algorithm \ref{alg:coeff_update4}                                    & Algorithm \ref{alg:coeff_update1}                                     & Algorithm \ref{alg:coeff_update2}                            & Algorithm \ref{alg:coeff_update4}                                                                                  \\
                &           & \multicolumn{2}{c}{Select best update every time (Algorithm \ref{alg:rID_overview})} & \frmt{24.31}                 & \frmt{19.64}                                 & \frmt{12.52}                                 & \frmt{9.35}                                  & \frmt{4.82}                                  & \frmt{3.94}                                                                               \\
        \bottomrule
\end{tabular}
        }
        \caption{The relative error of randomized ID with different coefficient computation algorithms and different numbers of new columns. The estimated error using single-pass Hutch++ is included in the bracket. The testing is performed on the ignition database.}
        \label{tab:ignition}
\end{table}

In Table \ref{tab:ignition}, we measure the effect of the target rank $k$ on the relative error of different coefficient computation algorithms and use the truncated SVD and single-pass randomized SVD \cite{yu2017single} as the benchmark for comparison. From the last row of each category in Table \ref{tab:ignition}, we still find that the reconstruction accuracy keeps improving when the target size increases. In terms of different coefficient algorithms, we find the reconstruction accuracy depends on the target rank setting and we cannot find an algorithm that beats the rest in every target rank. In addition, the Hutch++ error estimator still predicts the Frobenius error accurately and determines the coefficient matrix that produces a smaller reconstruction error. 

Regarding the performance of the ridge leverage score-based CSS algorithm, we first compare it with other CSS methods in a multi-pass setting. We observe similar results as in the channel flow dataset in that the CSS based on ridge leverage score (Algorithm \ref{alg:rID_overview}) performs slightly worse than the standard leverage score CSS and both methods perform slightly worse than the CPQR method. However, both leverage score-based methods perform significantly better than the residual-based CSS (Algorithm \ref{alg:CSS_residual}). 
A possible reason for the poor performance of the residual-based CSS (Algorithm \ref{alg:CSS_residual}) is that it cannot prune the already selected column basis and therefore fails to preserve the key time frames in the later stages of the dataset. 
Instead, the ridge leverage score-based method (Algorithm \ref{alg:rID_overview}) performs sampling among the entire data matrix and thus can effectively handle the dataset with a large number of time steps.

We then compare the ridge leverage score-based CSS (Algorithm \ref{alg:rID_overview}) with the residual-based CSS (Algorithm \ref{alg:CSS_residual}) in an online setting where the coefficient matrix is computed simultaneously with CSS. We observe that the ridge leverage score-based method outperforms the residual-based CSS method, which verifies the ability of our framework to handle datasets with drastic temporal variations.

Some other general observations are that the randomized SVD is now considerably sub-optimal compared to the usual SVD, whereas there was little difference in the turbulent flow data. Similarly to the turbulent flow data, the best ID method (the CPQR) is comparable to the SVD at low ranks (with a factor of $2\times$)) but much worse at high ranks (within a factor of $5\times$). Also as before, the best one-pass ID methods are comparable to CPQR at low rank, but a factor of $3\times$ or worse at high rank.  Hence, again we conclude that the online ID is within a small factor of optimal (the 2-pass SVD) at low ranks, but has more room for improvement at high ranks.

\FloatBarrier
\subsection{The NSTX Gas Puff Image (GPI) data}
\begin{figure}[!ht]
        \centering
        \includegraphics[width =\textwidth]{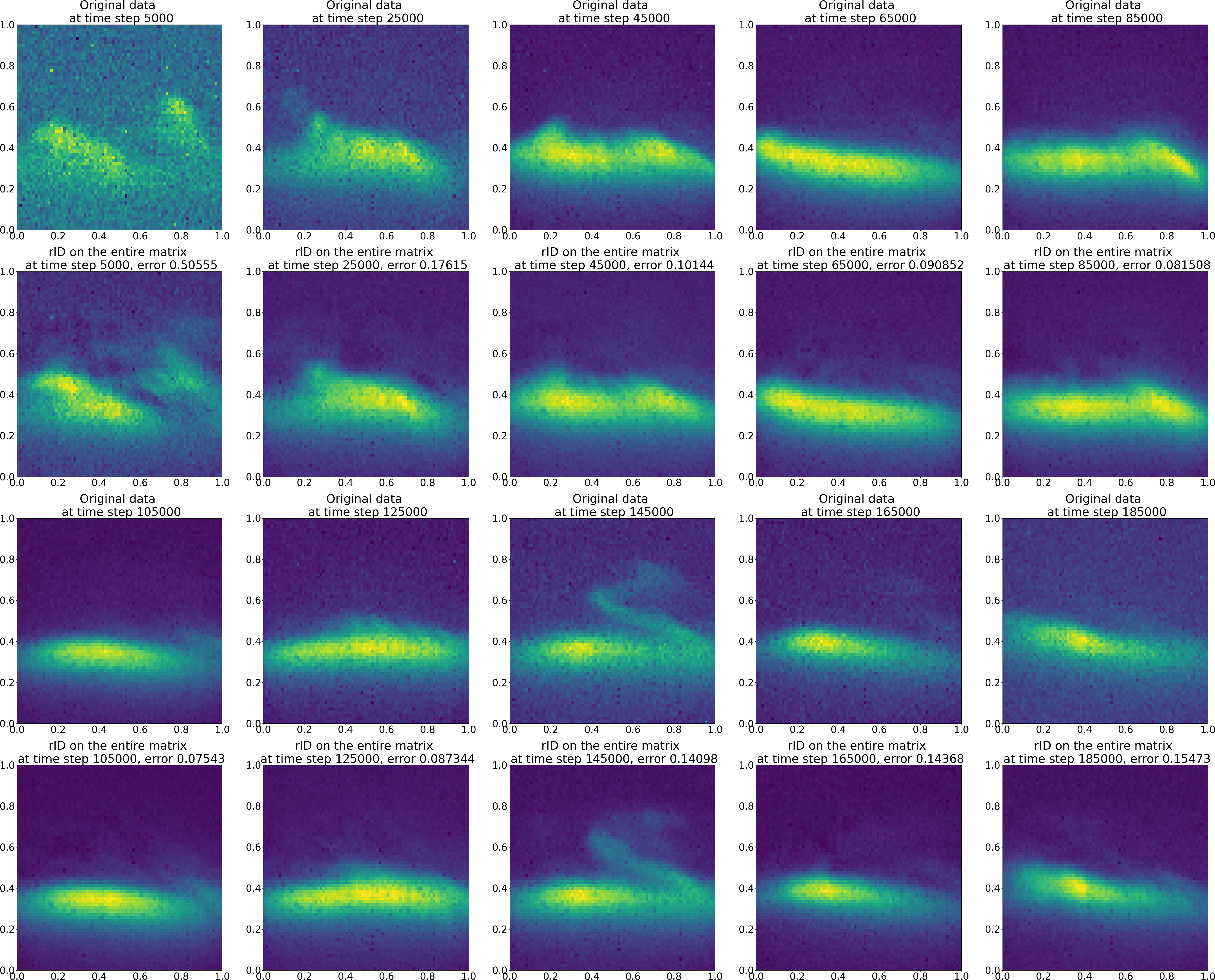}
        \caption{Rank $k=100$ reconstruction of NSTX gas puff image data over a grid of size $80 \times 64$ at different time instances. For each two rows, the top one represents the original data while the bottom one represents the reconstruction data.}
        \label{fig:nstxgpi}
\end{figure}

In the previous two tests, the input data matrix $m\times n$ data matrix $\bA$ satisfies $m \gg n$. To evaluate the performance of our online framework in handling the datasets with many time steps ($n\gg m$), we next consider the NSTX GPI dataset \cite{zhao2020sdrbench, SDRBench_git} that contains 369,357 data snapshots over a grid of size $80\times 64$. Since the dataset was obtained from an experiment, the first half of the time points are used to calibrate the measurement noise. Here, we compute the average of the first 169,357 frames as the baseline noise which we subtract from the entire dataset. Then, we use the remaining 200,000 frames as the input data to our framework. The reconstructed data with target rank $k=100$ at different time steps is shown in Fig.\ \ref{fig:nstxgpi}. We find that the reconstruction result effectively preserves the profile of the gas puff. In particular, we find that the original data at the initial time steps (5,000 and 25,000) is still noisy and the gas puff profile is hard to identify. However, the reconstruction manages to reduce the noise signal, demonstrating that our randomized ID method captures the key time frames as column bases for compression.
\begin{table}[!htb]
        \resizebox{\textwidth}{!}{
\centering
\begin{tabular}{ccllllll}
        \toprule
                          &           &                                                   &                              & \multicolumn{4}{c}{Target rank}                                                                                                                                                                  \\
        \cmidrule(l){5-8}
        Dataset           & \# passes & Column selection                                  & Coeff.\ update               & 100                                          & 200                                          & 400                                          & 800                                          \\
        \midrule
        \multirow{14}{*}{\shortstack{NSTX GPI data                                                                                                                                                                                                                                                  \\ $80 \times 64 \times 50,000$}}
                          & $2$       & \multicolumn{2}{c}{Truncated SVD}                 & \frmt{8.78}                  & \frmt{8.17}                                  & \frmt{7.61}                                  & \frmt{6.95}                                                                                 \\
                          & $1$       & \multicolumn{2}{c}{Randomized SVD}                & \frmt{10.12}                 & \frmt{9.35}                                  & \frmt{7.91}                                  & \frmt{7.35}                                                                                 \\
        \cmidrule(l){2-8}
                          & $k$       & CPQR                                        & Least squares (Eq.~\eqref{eq:coeff_ID2})  & \frmt{10.52}                                 & \frmt{9.75}                                  & \frmt{7.85}                                  & \frmt{6.54}                                  \\
                          & $2$       & Leverage score                                    & Least squares (Eq.~\eqref{eq:coeff_ID2})  & \frmt{13.17}                                 & \frmt{11.42}                                 & \frmt{8.96}                                  & \frmt{7.97}                                  \\
                          & $2$       & Residual-based (Algorithm \ref{alg:CSS_residual})                                & Least squares (Eq.~\eqref{eq:coeff_ID2})  & \frmt{36.34}                                 & \frmt{30.56}                                 & \frmt{19.16}                                 & \frmt{12.31}                                 \\
                          & $2$       & Ridge leverage score-based (Algorithm \ref{alg:rID_overview})                     & Least squares (Eq.~\eqref{eq:coeff_ID2})  & \frmt{11.25}                                 & \frmt{9.54}                                  & \frmt{9.12}                                  & \frmt{8.54}                                  \\
        \cmidrule(l){2-8}
                          & $1$       & Residual-based (Algorithm \ref{alg:CSS_residual})                                & Algorithm \ref{alg:coeff_update1}              & \frmt{39.71}(\frmt{40.16})                   & \best{\frmt{32.54}}(\frmt{34.83})          & \best{\frmt{22.63}}(\best{\frmt{21.53}}) & \frmt{16.12}(\frmt{19.62})                   \\
                          & $1$       & Residual-based (Algorithm \ref{alg:CSS_residual})                                & Algorithm \ref{alg:coeff_update2} & \frmt{38.96}(\frmt{39.81})                   & \frmt{33.78}(\best{\frmt{31.58}})          & \frmt{25.13}(\frmt{27.24})                   & \best{\frmt{20.83}}(\best{\frmt{17.34}}) \\
                          & $1$       & Residual-based (Algorithm \ref{alg:CSS_residual})                                & Algorithm \ref{alg:coeff_update3}     & \frmt{38.87}(\frmt{37.65})                   & \frmt{34.62}(\frmt{34.67})                   & \frmt{22.95}(\frmt{23.75})                   & \frmt{16.14}(\frmt{21.87})                   \\
                          & $1$       & Residual-based (Algorithm \ref{alg:CSS_residual})                                & Algorithm \ref{alg:coeff_update4}                    & \best{\frmt{37.79}}(\best{\frmt{36.10}}) & \frmt{33.02}(\frmt{32.17})                   & \frmt{24.24}(\frmt{23.64})                   & \frmt{16.05}(\frmt{20.34})                   \\
        \cmidrule(l){2-8}
                          & $1$       & Ridge leverage score-based (Algorithm \ref{alg:rID_overview})                     & Algorithm \ref{alg:coeff_update1}              & \frmt{13.12}(\frmt{14.21})                   & \frmt{11.24}(\frmt{11.39})                   & \frmt{9.36}(\frmt{9.53})                     & \best{\frmt{8.51}}(\best{\frmt{8.12}})   \\
                          & $1$       & Ridge leverage score-based (Algorithm \ref{alg:rID_overview})                     & Algorithm \ref{alg:coeff_update2} & \best{\frmt{12.91}}(\best{\frmt{13.21}}) & \frmt{11.43}(\frmt{11.52})                   & \frmt{9.57}(\frmt{10.27})                    & \frmt{8.92}(\frmt{9.25})                     \\
                          & $1$       & Ridge leverage score-based (Algorithm \ref{alg:rID_overview})                     & Algorithm \ref{alg:coeff_update3}     & \frmt{13.25}(\frmt{13.47})                   & \frmt{11.49}(\frmt{11.58})                   & \frmt{9.81}(\frmt{10.46})                    & \frmt{8.64}(\frmt{8.37})                     \\
                          & $1$       & Ridge leverage score-based (Algorithm \ref{alg:rID_overview})                     & Algorithm \ref{alg:coeff_update4}                    & \frmt{13.17}(\frmt{13.59})                   & \best{\frmt{11.15}}(\best{\frmt{11.22}}) & \best{\frmt{9.28}}(\best{\frmt{9.31}})   & \frmt{8.77}(\frmt{8.39})                     \\
        \cmidrule(l){2-8}
                          &           & \multicolumn{2}{c}{Update decision}               & Algorithm \ref{alg:coeff_update2} & Algorithm \ref{alg:coeff_update4}                                    & Algorithm \ref{alg:coeff_update4}                                    & Algorithm \ref{alg:coeff_update1}                                                                             \\
                          &           & \multicolumn{2}{c}{Select best update every time (Algorithm \ref{alg:rID_overview})} & \frmt{12.38}                 & \frmt{10.45}                                 & \frmt{8.99}                                  & \frmt{7.99}                                                                                 \\
        \midrule
        \multirow{14}{*}{\shortstack{NSTX GPI data                                                                                                                                                                                                                                                  \\ $80 \times 64 \times 100,000$}}
                          & $2$       & \multicolumn{2}{c}{Truncated SVD}                 & \frmt{11.33}                 & \frmt{10.66}                                 & \frmt{9.99}                                  & \frmt{9.18}                                                                                 \\
                          & $1$       & \multicolumn{2}{c}{Randomized SVD}                & \frmt{12.84}                 & \frmt{11.74}                                 & \frmt{10.65}                                 & \frmt{10.12}                                                                                \\
        \cmidrule(l){2-8}
                          & $k$       & CPQR                                        & Least squares (Eq.~\eqref{eq:coeff_ID2})  & \frmt{12.92}                                 & \frmt{11.35}                                 & \frmt{10.52}                                 & \frmt{9.03}                                  \\
                          & $2$       & Leverage score                                    & Least squares (Eq.~\eqref{eq:coeff_ID2})  & \frmt{13.74}                                 & \frmt{13.63}                                 & \frmt{12.08}                                 & \frmt{9.37}                                  \\
                          & $2$       & Residual-based (Algorithm \ref{alg:CSS_residual})                                & Least squares (Eq.~\eqref{eq:coeff_ID2})  & \frmt{37.43}                                 & \frmt{31.74}                                 & \frmt{23.65}                                 & \frmt{14.91}                                 \\
                          & $2$       & Ridge leverage score-based (Algorithm \ref{alg:rID_overview})                     & Least squares (Eq.~\eqref{eq:coeff_ID2})  & \frmt{13.64}                                 & \frmt{11.93}                                 & \frmt{11.21}                                 & \frmt{9.89}                                  \\
        \cmidrule(l){2-8}
                          & $1$       & Residual-based (Algorithm \ref{alg:CSS_residual})                                & Algorithm \ref{alg:coeff_update1}              & \frmt{38.74}(\frmt{39.62})                   & \frmt{32.50}(\frmt{34.13})                   & \best{\frmt{26.63}}(\best{\frmt{26.81}}) & \frmt{19.12}(\frmt{20.34})                   \\
                          & $1$       & Residual-based (Algorithm \ref{alg:CSS_residual})                                & Algorithm \ref{alg:coeff_update2} & \frmt{36.66}(\frmt{35.18})                   & \best{\frmt{31.54}}(\best{\frmt{31.29}}) & \frmt{29.13}(\frmt{29.42})                   & \best{\frmt{23.83}}(\best{\frmt{21.75}}) \\
                          & $1$       & Residual-based (Algorithm \ref{alg:CSS_residual})                                & Algorithm \ref{alg:coeff_update3}     & \frmt{37.25}(\frmt{37.24})                   & \frmt{33.21}(\frmt{32.70})                   & \frmt{26.95}(\frmt{27.98})                   & \frmt{20.42}(\frmt{24.78})                   \\
                          & $1$       & Residual-based (Algorithm \ref{alg:CSS_residual})                                & Algorithm \ref{alg:coeff_update4}                    & \best{\frmt{35.79}}(\best{\frmt{34.59}}) & \frmt{32.40}(\frmt{33.33})                   & \frmt{28.24}(\frmt{28.90})                   & \frmt{20.30}(\frmt{22.34})                   \\
        \cmidrule(l){2-8}
                          & $1$       & Ridge leverage score-based (Algorithm \ref{alg:rID_overview})                     & Algorithm \ref{alg:coeff_update1}              & \frmt{15.15}(\frmt{15.37})                   & \frmt{14.50}(\frmt{16.13})                   & \frmt{13.13}(\frmt{13.97})                   & \frmt{11.12}(\frmt{11.96})                   \\
                          & $1$       & Ridge leverage score-based (Algorithm \ref{alg:rID_overview})                     & Algorithm \ref{alg:coeff_update2} & \frmt{16.32}(\frmt{16.21})                   & \best{\frmt{13.54}}(\best{\frmt{14.29}}) & \best{\frmt{11.63}}(\best{\frmt{11.81}}) & \best{\frmt{10.83}}(\best{\frmt{10.34}}) \\
                          & $1$       & Ridge leverage score-based (Algorithm \ref{alg:rID_overview})                     & Algorithm \ref{alg:coeff_update3}     & \frmt{15.58}(\frmt{14.90})                   & \frmt{14.21}(\frmt{16.70})                   & \frmt{11.95}(\frmt{12.98})                   & \frmt{11.02}(\frmt{10.78})                   \\
                          & $1$       & Ridge leverage score-based (Algorithm \ref{alg:rID_overview})                     & Algorithm \ref{alg:coeff_update4}                    & \best{\frmt{14.79}}(\best{\frmt{14.30}}) & \frmt{15.40}(\frmt{15.33})                   & \frmt{12.24}(\frmt{13.12})                   & \frmt{11.30}(\frmt{11.54})                   \\
        \cmidrule(l){2-8} &           & \multicolumn{2}{c}{Update decision}               & Algorithm \ref{alg:coeff_update4}                    & Algorithm \ref{alg:coeff_update2}                 & Algorithm \ref{alg:coeff_update2}                 & Algorithm \ref{alg:coeff_update2}                                                                \\
                          &           & \multicolumn{2}{c}{Select best update every time (Algorithm \ref{alg:rID_overview})} & \frmt{14.47}                 & \frmt{12.77}                                 & \frmt{11.48}                                 & \frmt{10.45}                                                                                \\
        \midrule
        \multirow{14}{*}{\shortstack{NSTX GPI data                                                                                                                                                                                                                                                  \\ $80 \times 64 \times 150,000$}}
                          & $2$       & \multicolumn{2}{c}{Truncated SVD}                 & \frmt{13.36}                 & \frmt{12.64}                                 & \frmt{11.88}                                 & \frmt{10.94}                                                                                \\
                          & $1$       & \multicolumn{2}{c}{Randomized SVD}                & \frmt{14.31}                 & \frmt{13.24}                                 & \frmt{12.16}                                 & \frmt{11.75}                                                                                \\
        \cmidrule(l){2-8}
                          & $k$       & CPQR                                        & Least squares (Eq.~\eqref{eq:coeff_ID2})  & \frmt{14.95}                                 & \frmt{13.11}                                 & \frmt{12.38}                                 & \frmt{11.37}                                 \\
                          & $2$       & Leverage score                                    & Least squares (Eq.~\eqref{eq:coeff_ID2})  & \frmt{17.58}                                 & \frmt{15.62}                                 & \frmt{14.69}                                 & \frmt{12.37}                                 \\
                          & $2$       & Residual-based (Algorithm \ref{alg:CSS_residual})                                & Least squares (Eq.~\eqref{eq:coeff_ID2})  & \frmt{39.76}                                 & \frmt{32.39}                                 & \frmt{25.69}                                 & \frmt{17.14}                                 \\
                          & $2$       & Ridge leverage score-based (Algorithm \ref{alg:rID_overview})                     & Least squares (Eq.~\eqref{eq:coeff_ID2})  & \frmt{15.12}                                 & \frmt{13.95}                                 & \frmt{12.54}                                 & \frmt{11.96}                                 \\
        \cmidrule(l){2-8}
                          & $1$       & Residual-based (Algorithm \ref{alg:CSS_residual})                                & Algorithm \ref{alg:coeff_update1}              & \frmt{42.12}(\frmt{43.13})                   & \frmt{37.24}(\frmt{36.29})                   & \frmt{31.16}(\frmt{30.90})                   & \frmt{18.29}(\frmt{19.16})                   \\
                          & $1$       & Residual-based (Algorithm \ref{alg:CSS_residual})                                & Algorithm \ref{alg:coeff_update2} & \frmt{43.62}(\frmt{43.23})                   & \best{\frmt{36.10}}(\best{\frmt{37.12}}) & \frmt{30.31}(\frmt{30.42})                   & \best{\frmt{19.75}}(\best{\frmt{20.19}}) \\
                          & $1$       & Residual-based (Algorithm \ref{alg:CSS_residual})                                & Algorithm \ref{alg:coeff_update3}     & \frmt{42.25}(\frmt{42.57})                   & \frmt{38.21}(\frmt{38.16})                   & \frmt{29.75}(\frmt{31.98})                   & \frmt{21.62}(\frmt{22.91})                   \\
                          & $1$       & Residual-based (Algorithm \ref{alg:CSS_residual})                                & Algorithm \ref{alg:coeff_update4}                    & \best{\frmt{41.79}}(\best{\frmt{42.01}}) & \frmt{39.41}(\frmt{37.73})                   & \best{\frmt{29.63}}(\best{\frmt{29.81}}) & \frmt{21.83}(\frmt{21.45})                   \\
        \cmidrule(l){2-8}
                          & $1$       & Ridge leverage score-based (Algorithm \ref{alg:rID_overview})                     & Algorithm \ref{alg:coeff_update1}              & \frmt{17.17}(\frmt{18.30})                   & \frmt{15.95}(\frmt{15.12})                   & \frmt{14.05}(\frmt{14.30})                   & \frmt{13.39}(\frmt{14.81})                   \\
                          & $1$       & Ridge leverage score-based (Algorithm \ref{alg:rID_overview})                     & Algorithm \ref{alg:coeff_update2} & \frmt{17.26}(\frmt{18.14})                   & \best{\frmt{14.18}}(\best{\frmt{14.30}}) & \frmt{14.19}(\frmt{14.44})                   & \frmt{13.17}(\frmt{14.07})                   \\
                          & $1$       & Ridge leverage score-based (Algorithm \ref{alg:rID_overview})                     & Algorithm \ref{alg:coeff_update3}     & \frmt{17.38}(\frmt{17.72})                   & \frmt{14.92}(\frmt{14.74})                   & \best{\frmt{13.76}}(\best{\frmt{13.92}}) & \frmt{13.24}(\frmt{13.64})                   \\
                          & $1$       & Ridge leverage score-based (Algorithm \ref{alg:rID_overview})                     & Algorithm \ref{alg:coeff_update4}                    & \best{\frmt{17.04}}(\best{\frmt{17.17}}) & \frmt{14.67}(\frmt{14.83})                   & \frmt{14.85}(\frmt{14.25})                   & \best{\frmt{13.13}}(\best{\frmt{12.96}}) \\
        \cmidrule(l){2-8}
                          &           & \multicolumn{2}{c}{Update decision}               & Algorithm \ref{alg:coeff_update4}                    & Algorithm \ref{alg:coeff_update2}                 & Algorithm \ref{alg:coeff_update3}                     & Algorithm \ref{alg:coeff_update4}                                                                                   \\
                          &           & \multicolumn{2}{c}{Select best update every time (Algorithm \ref{alg:rID_overview})} & \frmt{16.28}                 & \frmt{14.65}                                 & \frmt{13.45}                                 & \frmt{12.37}                                                                                \\
        \midrule
        \multirow{14}{*}{\shortstack{NSTX GPI data                                                                                                                                                                                                                                                  \\$80 \times 64 \times 200,000$}}
                          & $2$       & \multicolumn{2}{c}{Truncated SVD}                 & \frmt{14.98}                 & \frmt{14.19}                                 & \frmt{13.34}                                 & \frmt{12.29}                                                                                \\
                          & $1$       & \multicolumn{2}{c}{Randomized SVD}                & \frmt{16.35}                 & \frmt{15.06}                                 & \frmt{14.35}                                 & \frmt{13.15}                                                                                \\
        \cmidrule(l){2-8}
                          & $k$       & CPQR                                        & Least squares (Eq.~\eqref{eq:coeff_ID2})  & \frmt{14.15}                                 & \frmt{15.13}                                 & \frmt{14.24}                                 & \frmt{12.77}                                 \\
                          & $2$       & Leverage score                                    & Least squares (Eq.~\eqref{eq:coeff_ID2})  & \frmt{17.26}                                 & \frmt{17.32}                                 & \frmt{15.25}                                 & \frmt{12.11}                                 \\
                          & $2$       & Residual-based (Algorithm \ref{alg:CSS_residual})                                & Least squares (Eq.~\eqref{eq:coeff_ID2})  & \frmt{40.59}                                 & \frmt{33.83}                                 & \frmt{26.02}                                 & \frmt{18.35}                                 \\
                          & $2$       & Ridge leverage score-based (Algorithm \ref{alg:rID_overview})                     & Least squares (Eq.~\eqref{eq:coeff_ID2})  & \frmt{16.90}                                 & \frmt{15.05}                                 & \frmt{14.51}                                 & \frmt{13.16}                                 \\
        \cmidrule(l){2-8}
                          & $1$       & Residual-based (Algorithm \ref{alg:CSS_residual})                                & Algorithm \ref{alg:coeff_update1}              & \best{\frmt{39.86}}(\best{\frmt{40.11}}) & \frmt{34.97}(\frmt{35.62})                   & \best{\frmt{26.11}}(\best{\frmt{27.94}}) & \frmt{19.95}(\frmt{21.76})                   \\
                          & $1$       & Residual-based (Algorithm \ref{alg:CSS_residual})                                & Algorithm \ref{alg:coeff_update2} & \frmt{43.62}(\frmt{44.78})                   & \best{\frmt{34.10}}(\best{\frmt{33.29}}) & \frmt{27.31}(\frmt{29.64})                   & \frmt{20.83}(\frmt{21.35})                   \\
                          & $1$       & Residual-based (Algorithm \ref{alg:CSS_residual})                                & Algorithm \ref{alg:coeff_update3}     & \frmt{41.25}(\frmt{42.49})                   & \frmt{34.66}(\frmt{34.57})                   & \frmt{26.75}(\frmt{27.83})                   & \frmt{20.62}(\frmt{22.17})                   \\
                          & $1$       & Residual-based (Algorithm \ref{alg:CSS_residual})                                & Algorithm \ref{alg:coeff_update4}                    & \frmt{42.71}(\frmt{43.26})                   & \frmt{35.91}(\frmt{35.36})                   & \frmt{28.85}(\frmt{29.97})                   & \best{\frmt{19.07}}(\best{\frmt{20.38}}) \\
        \cmidrule(l){2-8}
                          & $1$       & Ridge leverage score-based (Algorithm \ref{alg:rID_overview})                     & Algorithm \ref{alg:coeff_update1}              & \frmt{18.33}(\frmt{18.25})                   & \frmt{17.32}(\frmt{18.74})                   & \frmt{15.23}(\frmt{16.39})                   & \frmt{14.94}(\frmt{14.81})                   \\
                          & $1$       & Ridge leverage score-based (Algorithm \ref{alg:rID_overview})                     & Algorithm \ref{alg:coeff_update2} & \frmt{18.94}(\frmt{18.13})                   & \best{\frmt{17.34}}(\best{\frmt{17.56}}) & \frmt{15.84}(\frmt{16.13})                   & \frmt{14.73}(\frmt{15.17})                   \\
                          & $1$       & Ridge leverage score-based (Algorithm \ref{alg:rID_overview})                     & Algorithm \ref{alg:coeff_update3}     & \frmt{18.81}(\frmt{18.24})                   & \frmt{17.23}(\frmt{17.91})                   & \best{\frmt{15.13}}(\best{\frmt{15.96}}) & \frmt{14.85}(\frmt{15.29})                   \\
                          & $1$       & Ridge leverage score-based (Algorithm \ref{alg:rID_overview})                     & Algorithm \ref{alg:coeff_update4}                    & \best{\frmt{18.24}}(\best{\frmt{17.87}}) & \frmt{17.16}(\frmt{17.93})                   & \frmt{15.39}(\frmt{16.29})                   & \best{\frmt{14.32}}(\best{\frmt{14.52}}) \\
        \cmidrule(l){2-8}
                          &           & \multicolumn{2}{c}{Update decision}               & Algorithm \ref{alg:coeff_update4}                    & Algorithm \ref{alg:coeff_update2}                 & Algorithm \ref{alg:coeff_update3}                     & Algorithm \ref{alg:coeff_update4}                                                                                   \\
                          &           & \multicolumn{2}{c}{Select best update every time (Algorithm \ref{alg:rID_overview})} & \frmt{18.00}                 & \frmt{16.95}                                 & \frmt{15.11}                                 & \frmt{13.93}                                                                                \\
        \bottomrule
\end{tabular}
        }
        \caption{The relative error of randomized ID with different coefficient computation algorithms and different numbers of new columns. The estimated error using single-pass Hutch++ is included in the bracket. Testing is performed on NSTX Gas Puff Image Data after denoising.}
        \label{tab:gpidata}
\end{table}

In Table \ref{tab:gpidata}, we monitor the reconstruction accuracy during the streaming data process and summarize the accuracy at four different time steps: 50,000, 100,000, 150,000, and 200,000. For each case, we measure the effect of the target rank ($k$) coefficient computation algorithm on the relative error and use the SVD and a single-pass randomized SVD as the benchmarks for comparison. Similar to the previous results, the reconstruction accuracy keeps improving as the target rank increases. In addition, no single coefficient calculation strategy leads to more accurate reconstructions uniformly across the ranks, and the NA-Hutch++ scheme produces accurate error estimates. By selecting the best updating algorithm via the NA-Hutch++ error estimator, whenever the column basis is updated, we obtain more accurate reconstructions than using a single coefficient computation algorithm.

Regarding the performance of the ridge leverage score-based CSS algorithm, we first compare it with other CSS methods in a two-pass setting. We observe similar results as in the channel flow dataset where the CSS method performs slightly worse than the standard leverage score CSS \cite{jolliffe1972discarding} and both leverage score-based CSS methods perform slightly worse than the CPQR method. And as before, both leverage score-based methods perform significantly better than the residual-based CSS. By comparing our choice of CSS method with the residual-based CSS in a single-pass setting, we similarly observe that the ridge leverage score-based method outperforms the residual-based CSS method, which demonstrates the ability of our approach to handle a dataset consisting of many time frames.

Unlike the previous two datasets, the relative performance of rID in comparison to the optimal (2 pass SVD) does not change much with rank, and is always almost comparable. Notably, the SVD error at high rank is much worse than the previous datasets, possibly due to the noise. This suggests that the ID, and specifically our online rID, is (relatively) more robust to noise than the SVD, and particularly useful in high-noise, low-accuracy regimes.


\section{Conclusions}
\label{sec:Conclusion}
In this study, we develop a randomized ID approach to perform temporal compression in a single pass over large-scale, streaming simulation data. 
We demonstrate that the proposed algorithms can be employed to obtain the low-rank approximation of data taken from (1) a 3D turbulent channel flow; (2) a 2D ignition problem; and (3) the NSTX Gas Puff Image (GPI) measurements. We also present an adaptation of the single-pass, randomized ID method that improves the reconstruction error of the solution gradient, when desired, without losing considerable accuracy on the reconstructed solution itself. The numerical results show clearly that of the two reasonable streaming strategies for picking the column basis, the ridge leverage score-based method we use in our main algorithm  (Algorithm \ref{alg:rID_overview}) is consistently much better than the residual-based method (Algorithm \ref{alg:CSS_residual}). As for the type of sketching (Algorithms \ref{alg:coeff_update1},\ref{alg:coeff_update2},\ref{alg:coeff_update3},\ref{alg:coeff_update4}) for the coefficient update, there was not one method that always dominated the others. Since our error estimate is reasonably accurate, it allowed us to determine on-the-fly which strategy was most accurate and use it.

Comparing our streaming ID to either an offline ID or an offline SVD, the first two datasets showed that the streaming ID was nearly as good as the offline ID (which in turn was nearly as good as the SVD) for low ranks, but not for high ranks. In contrast, the noisier third dataset, which had much larger $n$ and hence all the ranks were considerably larger, showed that the rID was nearly as good as the SVD throughout all ranks studied. We conclude that for datasets of this last type, the streaming ID is often more attractive than the SVD since it has similar error yet maintains interpretability. We also see that for datasets that require further derived quantities (such as deriving the vorticity from the gradient of the velocity field), explicitly taking this into account using discrete gradients can improve not just the quality of the derived quantity but also, surprisingly, the original field as well.

Future extensions of the proposed strategy include parallel, \textit{in situ} implementation in conjunction with the PDE solver generating solution data. In addition, the current rank-based implementation can be altered to an error-based formulation where the solution rank is increased to achieve a prescribed accuracy, which is attainable. Finally, the proposed temporal compression can be augmented with spatial compression by utilizing spatial compressors, such as SZ \cite{di2016fast} and FPZIP \cite{lindstrom2006fast}, to reduce the column basis vectors, thereby achieving higher compression factors. This can be done \emph{before} coefficient updates so that slightly adjusted coefficients can partially account for the reduced accuracy in basis vectors.

\section*{Acknowledgements}
This material is based upon work supported by the Department of Energy Advanced Scientific Computing Research Award DE-SC0022283. AD's work has also been partially supported by AFOSR grant FA9550-20-1-0138 and SB's work has been partially supported by DOE award DE-SC0023346. We would also like to thank Kenneth Jansen, John Evans, and Jeff Hadley from the University of Colorado Boulder for their helpful discussions surrounding this work.

\section*{Declaration of competing interest}
The authors declare that they have no known competing financial interests or personal relationships that could have appeared to influence the work reported in this paper.

\FloatBarrier
\bibliographystyle{elsarticle-num}
\bibliography{sample}

\end{document}